\documentclass[10pt,reqno]{article}
\usepackage{amssymb, amsmath, amsthm, amsfonts, amscd, epsfig}
\usepackage{bbm}
\usepackage{placeins}
\usepackage[english]{babel}
\usepackage{bm}
\usepackage{graphicx}
\usepackage{setspace}

\newtheorem{theorem}{Theorem}[section]
\newtheorem{cor}[theorem]{Corollary}
\newtheorem{lemma}[theorem]{Lemma}
\newtheorem{prop}[theorem]{Proposition}

\newcommand{\nm}{\noalign{\smallskip}}

\def\ep{\epsilon}

\newcommand{\Bp}{\mathbf{p}}

\newcommand{\Bc}{\mathbf{c}}

\newcommand{\Bx}{\mathbf{x}}

\newcommand{\RR}{\mathbb{R}}
\newcommand{\CC}{\mathbb{C}}
\newcommand{\NN}{\mathbb{N}}

\newcommand{\Ical}{\mathcal{I}}
\newcommand{\p}{\partial}

\newcommand{\pd}[2]{\frac {\p #1}{\p #2}}
\newcommand{\ds}{\displaystyle}
\newcommand{\eqnref}[1]{(\ref {#1})}
\newcommand{\pf}{\medskip \noindent {\sl Proof}. \ }
\renewcommand{\qed}{\hfill $\Box$ \medskip}

\newcommand{\beq}{\begin{equation}}
\newcommand{\eeq}{\end{equation}}
\newcommand{\Br}{\mathbf{r}}

\newcommand{\rone}{\tilde{r}_1}
\newcommand{\rtwo}{\tilde{r}_2}

\numberwithin{equation}{section}
\numberwithin{figure}{section}

\begin{document}

\title{Asymptotic analysis for superfocusing of the electric field in between two nearly touching metallic spheres\thanks{\footnotesize
This work is supported by the Korean Ministry of Science, ICT and Future Planning through NRF grant No. NRF-2013R1A1A3012931 (to M.L) and by the Korean Ministry of Education, Sciences and Technology through NRF grant No. 2012003224 (to S.Y)}}

\author{Mikyoung
Lim\thanks{\footnotesize Department of Mathematical Sciences,
Korea Advanced Institute of Science and Technology, Daejeon
305-701, Korea (mklim@kaist.ac.kr, shyu@kaist.ac.kr).} \and
Sanghyeon Yu\footnotemark[2]}

\maketitle
\begin{abstract}
We consider the enhancement of electric field in the presence of two perfectly conducting spheres. When the two spheres get closer, the electric field have a much larger magnitude compared to the external field in the small gap region between the two spheres. The enhanced field can be arbitrary large with the generic blow-up rate $|\ep\ln\ep|^{-1}$ in three dimensional space, where $\ep$ is the distance between the spheres.
In this paper we derive rigorously an asymptotic formula of the electric field consisting of elementary functions. The asymptotic formula explicitly characterizes superfocusing of the electric field in terms of the spheres radii, the distance between the spheres, and the external field. We illustrate our results with numerical calculations.
\end{abstract}



\section{Introduction}
Two nearly touching metallic spheres cause the enhancement of the electric field. In an external electric field of long wavelength compared to the size of spheres, the presence of nearly touching metallic spheres induces a very large electric field confined in the narrow gap region between the spheres. Since the field is concentrated in a small region compared to the wavelength of external field, this effect is often called the superfocusing. The superfocusing in nearly touching metallic spheres has attracted considerable attention due to its application to various imaging modalities such as the surface-enhanced raman spectroscopy (SERS) and the single molecule detection \cite{Pen, Rom}.

In this paper we formulate and analyze the superfocusing of the electric field in between two nearly touching metallic spheres. When a extremely low-frequency eternal field is applied, metallic objects behave like perfect conductors according to the Drude model for metals and it is valid to consider the quasi-static approximation, {\it i.e.} the Laplace's equation for the electric potential.
We assume that $B_1$ and $B_2$ are two perfectly conducting spheres embedded in $\RR^3$, which is occupied by the homogeneous material of the conductivity 1. Then we consider the electric potential $u$ which satisfies the following conductivity equation:
\beq\label{u:eqn}
\quad \left\{
\begin{array}{ll}
\ds\Delta u=0 \quad& \mbox{in } \mathbb{R}^3 \setminus \overline{B_1 \cup B_2}, \\
\ds u= \mbox{constant} \quad& \mbox{on }\p B_j, j=1,2,\\
\nm
\ds \int_{\p B_j} \partial_\nu u ~d\sigma=0,\ \quad& j=1,2,\\
\ds u(\Bx)-H(\Bx)=O(|\Bx|^{-2}) \quad&\mbox{as } |\Bx| \to \infty,
\end{array}
\right.
\end{equation}
where $H$ is the external electric potential which is given by an entire harmonic function.
Here and throughout $\nu$ and $\partial_\nu u$ respectively denote the outward unit normal vector to $\p B_j$ and the outward normal derivative of $u$ on $\p B_j$ ($j=1,2$).
The main goal of this paper is to understand rigorously the superfocusing of $\nabla u$ by deriving the asymptotic formula for $\nabla u$ as the distance $\ep$ between two spheres tends to 0.

The problem of the electrostatic interaction between two perfectly conducting spheres dates back to at least 1890s  when Maxwell calculated the electric potential energy of two charged conducting spheres \cite{Max}.
There are two classical methods to derive the exact solution for the electric potential. The first is the method of separation of variables in the bispherical coordinates \cite{Dav, Jef, Smy} and the second is the method of image charges \cite{Smy}.
Both methods express the electric potential as an infinite series which converges fast when the two spheres are well-separated.
However, in the case that the gap between the spheres are small, the solution series converges very slowly and it causes the difficulty in computing the solution accurately. Actually, the magnitude of the electric field may blow up to infinity as the distance $\ep$ tends to zero \cite{keller}.
The asymptotic behavior of the electric field in between two closely located conductors has been studied extensively in relation with the computation of the effective conductivity in composite materials \cite{Bat,keller,McP,Pol}.
In \cite{McP}, McPhedran and colleagues considered two nearly touching cylinders of highly conducting materials. There, they approximated discrete image charges by an continuous charge distribution. And, based on this approximation, they derived asymptotics for the multipole coefficients of the electric potential and computed the effective conductivity of the composite material which consists of densely packed arrays of highly conducting cylinders.
For three dimensional case, Poladian obtained similar result for highly conducting spheres in \cite{Pol,Pol2}. It is worth to mention that this method was extended to the two-dimensional linear elasticity \cite{McP2}.

Lately, it has been intensively studied the singular behavior of the electric field, which is the gradient of the solution $u$ to \eqnref{u:eqn} in this paper. It was shown that $|\nabla u|$ is bounded independently of $\ep$ when the conductivities of embedded inclusions are finite and strictly positive \cite{LN, LV}. However, if the conductivities of the inclusions degenerate to $\infty$ (perfectly conducting), then the gradient may blow up as $\ep$ tends to 0. The generic rate of the gradient blow-up is $|\ep\ln\ep|^{-1}$ in three dimensions \cite{BLY, BLY2,KLY2, LY}, while it is $\ep^{-1/2}$ in two dimensions \cite{AKLLL, AKLLZ, AKL,bab, BLY,BLY2,BC, keller, Y, Y2}. The insulating case has the same blow-up rate as the perfectly conducting case in two dimensions. The gradient may or may not blow up depending on the given entire harmonic function $H$. In two dimensions, it was shown in \cite{AKLLZ} that the gradient may blow up only when the linear term of $H$ is nonzero.

Let us fix some notations to state the related results in details.
Since the Laplace's equation is invariant under rotation and shifting, we can denote $B_1$ and $B_2$ as
\begin{align}\notag B_1&=B(\Bc_1,r_1),\quad B_2=B(\Bc_2,r_2),\quad\Bc_j=(0,0,c_j),\ j=1,2,\end{align}
where $$
c_1=\frac{r_2^2-r_1^2-(r_1+r_2+\ep)^2}{2(r_1+r_2+\ep)}\quad\mbox{and}\quad c_2=c_1+r_1+r_2+\ep.
$$
Here $B(\mathbf{c},r)$ means the ball centered at $\Bc$ with radius $r$. The radii $r_1$ and $r_2$ can be different from each other.
We let  $\Bp_1\in B_1$ and $\Bp_2\in B_2$ be, respectively, the fixed points of combined reflections $R_1\circ R_2$ and $R_2\circ R_1$, where $R_j$ is the reflection w.r.t. $\p B_j$, {\it i.e.},
\beq\notag R_j(\mathbf{x})=\frac{r_j^2(\Bx-{\Bc}_j)}{|\Bx-{\Bc}_j|^2}+{\Bc}_j,\ j=1,2.\eeq
It can be easily shown that \begin{equation}\label{assump}\mathbf{p}_1=(0,0,-a_\ep)\quad\mbox{and}\quad\mathbf{p}_2=(0,0,a_\ep),
\end{equation}
where
\beq\label{def:a}a_\ep=\frac{\sqrt\ep\sqrt{ (2 r_1 + \ep) (2 r_2 + \ep) (2 r_1 + 2 r_2 +
  \ep)}}{2 (r_1 + r_2 + \ep)}=\sqrt{\frac{2r_1r_2}{r_1+r_2}}\sqrt{\ep}+O(\ep\sqrt\ep).\end{equation}

The blow-up behavior of the electric field can be characterized by the solution to the following equation:
\beq\label{h:eqn}
\quad \left\{
\begin{array}{ll}
\ds\Delta h=0 \quad& \mbox{in } \mathbb{R}^3 \setminus \overline{B_1 \cup B_2}, \\
\ds h= \mbox{constant}\quad& \mbox{on }\p B_j,\ j=1,2,\\
\nm
\ds \int_{\p B_j} \partial_\nu h ~ds=(-1)^{j+1}\quad&\mbox{for } j=1,2,\\
\ds h(\Bx)=O(|\Bx|^{-2}) \quad&\mbox{as } |\Bx| \to \infty.
\end{array}
\right.
\end{equation}
We call $h$ the singular function to \eqnref{u:eqn}.
It was derived in \cite{KLY} that the solution $u$ to \eqnref{u:eqn} can be decomposed into the singular and regular parts:
\beq\label{bdd_g_def}
u(\Bx)=C_H^\ep h(\Bx) + H(\Bx)+ r(\Bx) \quad \mbox{with }C_H^\ep=\frac{u|_{\p B_1}-u|_{\p B_2}}{h|_{\p B_1}-h|_{\p B_2}},
\eeq
where $\|\nabla r\|_\infty$ is bounded independently of $\ep$. Throughout this paper, the symbol $\|\cdot\|_\infty$ denotes $\|\cdot\|_{L^\infty(\RR^3\setminus\overline{(B_1\cup B_2)})}$. Once $h$ is obtained, one can consequently compute the asymptotic of $\nabla u(\Bx)$ by differentiating the right-hand side in Eq.\;\eqnref{bdd_g_def} and that of $C_H^\ep$ by applying the following relation obtained in \cite{Y, Y2}:
\beq\label{0417:udiff}u|_{\p B_1}-u|_{\p B_2}=\int_{\p B_1\cup \p B_2}H\partial_\nu h~d\sigma.\eeq
If $B_1$ and $B_2$ are either disks in two dimensional space or balls in three dimensional space, then the singular function $h$ is a potential function generated by two (for disks) or a sequence (for balls) of point charges \cite{LY}. It is worth to mention that the decomposition \eqnref{bdd_g_def} and \eqnref{0417:udiff} holds for general shaped inclusions.  In \cite{ACKLY, KLeeY}, it was obtained the gradient blow-up term of $u$ in terms of the solution to \eqnref{h:eqn} corresponding to the disks osculating to $B_j$'s when $B_j$'s are of convex shape.

For spherical perfect conductors in $\RR^3$, $h$ has been expressed as the electric potential generated by a sequence of point charges located at multiply reflected points with respect to the two spheres and, based on this expansion, upper and lower bounds of $|\nabla u|$ were obtained \cite{LY}. It was further investigated in \cite{KLY2} to derive an asymptotic formula for the same radius $r_1=r_2=r$: for $\Bx =(x_1,x_2,x_3)$ outside the two spheres, the solution $u$ satisfies
\beq\label{KLY1_b}
\nabla u(\Bx)=\frac{{C^\ep_H}}{\pi|\ln\ep|(\ep +rx_1^2+rx_2^2)}(\mathbf{e}_3 +\eta(\Bx))+\nabla g\quad\mbox{if }|(x_1,x_2)|\leq\frac{r}{|\ln\ep|^2},\eeq
where $\mathbf{e}_3=(0,0,1)$, $\|\nabla g\|_\infty$ is bounded regardless of $\ep$ and $|\eta(\Bx)|=O(|\ln\ep|^{-1})$,
and the concentration factor $C_H^\ep$ satisfies \beq\label{KLY1_a}C_H^\ep= 2\pi\sum_{n=1}^\infty \frac{r}{n}\left(H\Big(0,0,\frac{r}{n}\Big)-H\Big(0,0,-\frac{r}{n}\Big)\right)+O(\sqrt{\ep}|\ln\ep|).\eeq

While the equation \eqnref{KLY1_b} provides an asymptotic of $\nabla u$, the blow-up phenomenon of the electric field requires further investigation in view of the fact that, firstly, the unidentified function $\eta(\Bx)$ and the remainder term in \eqnref{KLY1_a} can also cause the blow-up and, secondly, the valid region for \eqnref{KLY1_b} degenerates to a point as $\ep$ tends to 0. It is worth to remark that the formula \eqnref{KLY1_b} is slightly modified from that in \cite{KLY2} to be valid for any positive number $r$, not just 1.

 In this paper, we derive an asymptotic of $\nabla u$ which completely characterizes the blow-up of the electric field due to the presence of the two nearly touching metallic spheres accepting different radii. The main results are as follows:
 \begin{itemize}
  \item[(i)]
We show that the remainder term in Eq.\;\eqnref{KLY1_a} (and in the modified equation which is valid for spheres of different radii) is actually $O(\ep|\ln\ep|)$. As an immediate consequence, we can replace $C_H^\ep$ by its limit as $\ep$ tends to zero, say $C_H$, in the decomposition \eqnref{bdd_g_def} of $u$ and in the asymptotic of $\nabla u$.  Furthermore, we calculate the series summation of index $n$ in Eq.\;\eqnref{KLY1_a} and completely rewrite it as the summation in terms of homogeneous polynomial order of the external field $H$.  This reformulated expression  of $C_H$ gives directly the necessary and sufficient condition for the gradient blow-up occurrence in terms of the two spheres radii, the distance $\ep$ between two spheres, and the external field $H$, see Theorem \ref{lem:CH} and Corollary \ref{cor:blowup}.

\item[(ii)] We provide the asymptotic formula of $\nabla u$ which is valid in the whole exterior region of the two spheres. The blow-up term is expressed explicitly in terms of elementary functions with coefficients depending on the spheres radii, the distance between the spheres, and the external field, see Theorem \ref{main_thm2}.

\item[(iii)] We identify the location, size and shape of the region where the gradient blow-up occurs.
The dimension of the blow-up region turns out to be of order $|\ln\ep|^{-1/2}$ which shrinks to a point as $\ep$ tends to zero. Moreover, this order $|\ln\ep|^{-1/2}$ is shown to be optimal. In other words, we prove and characterize the occurrence of superfocusing of the electric field, see Theorem \ref{prop:superfocusing}.

\end{itemize}

The main ingredients in this paper are the bispherical coordinate system and the Euler-Maclaurin formula. In bispherical coordinates, the exact solution $h$ of the problem \eqnref{h:eqn} can be obtained by the method of separation of variables. By some manipulations, such as changing the order of summations, on the exact series solution, we obtain a new series of Riemann sum type. We then apply the Euler-Maclaurin formula to approximate this new series by an integral. Based on the integral expression, we investigate the blow-up feature of the electric field. This approach originally comes from the previous work by the authors \cite{LY3}. There, it was derived the asymptotics of $u$ when two circular cylinders with finite conductivities in $\mathbb{R}^2$ are closely located.
  It is worth to mention that the asymptotic of the potential difference between two spheres has been derived when an uniform external field is applied \cite{lekner11}.

From our analysis, it turns out that the potential $u$ can be approximated by the integral of a piecewise continuous image charge distribution. This charge distribution is of similar form to that obtained in \cite{Pol}
up to a multiplicative constant. There, the image charge distribution was assumed to be continuous motivated by the physical intuition and was determined to satisfy several functional equations derived from boundary conditions on the spheres and additional physical assumptions.
In this paper, however, the continuous image charge distribution is derived rigorously without any physical assumptions. We emphasize that, thanks to the analysis with mathematical rigor, we are able not only to approximate the image charge distribution but also to go further to extract the blow-up term of the electric field.

The paper is organized as follows. In section 2, we state the main results. Section 3 is to review the definition and properties of the bispherical coordinate system. We provide two different series expansion for $h$ by the bispherical coordinates and estimate the concentration factor $C_H^\ep$ in section 4. The asymptotic formula for $\nabla u$ is derived in section 5. Section 6 provides two applications of the Euler-Maclaurin formula, and we illustrate the main results with numerical calculation in section 7. The conclusion is provided in section 8.

\section{Main Results}
In this section we fix some notations and state the main results.

We first denote
\beq\label{def:tilder}\tilde{r}=\frac{r_1 r_2}{r_1+r_2},\quad\tilde{r}_j=\frac{r_j}{r_1+r_2},\  j=1,2,\eeq
and \beq
\begin{cases}
\ds\mu_\ep=\frac{1}{2\pi\tilde{r}}\left[ |\ln \ep|+\ln \tilde{r}+\ln 2-2\frac{\psi_0(\tilde{r}_1)\psi_0(\tilde{r}_2)-\gamma^2}{\psi_0(\tilde{r}_1)+\psi_0(\tilde{r}_2)+2\gamma}\right]^{-1},\label{mu1002}\\[3mm]
\ds\mu_j=\frac{\psi_0(\tilde{r}_j)+\gamma}{\psi_0(\tilde{r}_1)+\psi_0(\tilde{r}_2)+2\gamma},\quad j=1,2,
\end{cases}
\eeq
where $\psi_0$ is the digamma function and $\gamma$ is the Euler's constant. We then define for $k\in\NN$ that
\beq
\label{def:mu}
\ds
\mathcal{Q}_k(r_1,r_2):=
4\pi\tilde{r}^{k+1}\left[
\Bigr(\mu_1+(-1)^{k+1}\mu_2\Bigr)\zeta(k+1)+\frac{\mu_1 \psi_k(\rtwo)+(-1)^{k+1}\mu_2 \psi_k(\rone)}{{k!}}\right],\eeq
where
$\psi_k$ is the polygamma function of order $k$ and $\zeta$ is the Riemann zeta function,{\it i.e.},
 $\zeta(z)=\sum_{n=1}^\infty \frac{1}{n^z}$ for $z\in \CC,\ \Re(z)>1.$ Here the symbol $\Re(\cdot)$ means the real part of complex number.
  The polygamma function is defined as the derivatives of the Gamma function $\Gamma$,{\it i.e.},
$$\psi_k(z)=\frac{d^{k+1}}{dz^{k+1}}\ln \Gamma(z),\quad z\in\CC,\ k=0,1,\dots.$$
 For $z\neq 0,-1,-2,\dots$, the polygamma function can be expressed as (see \cite{AS})
\beq\label{polygamma}
\psi_k(z)=
\ds\begin{cases}\ds(-1)^{k+1} k! \sum_{m=0}^\infty {(m+z)^{-k-1}}\quad&\mbox{for }k\geq 1\\
\ds-\gamma+\sum_{k=0}^\infty \frac{z-1}{(k+1)(k+z)}\quad&\mbox{for }k=0.
\end{cases}
\eeq

 The first main result is the limiting behavior of $C_H^\ep$. We give the proof of Theorem \ref{lem:CH} in section \ref{section:CH}.
\begin{theorem}\label{lem:CH}
\begin{itemize}
\item[\rm(a)]
The concentration factor $C_H^\ep$ satisfies
\beq \label{eqn:CH}C_H^\ep =C_H+O(\ep|\ln\ep|)\quad\mbox{as }\ep\mbox{ tends to zero},\eeq
where
\begin{align}
\ds {C}_H&=\
4\pi\mu_1\sum_{m=0}^\infty \bigg[ \frac{\tilde r}{m+1}H\Bigr(0,0,\frac{\tilde r}{m+1}\Bigr)
- \frac{\tilde r}{m+\tilde r_2}H\Bigr(0,0,\frac{-\tilde r}{m+\tilde r_2}\Bigr)\bigg]\notag\\
&\ds
+4\pi\mu_2\sum_{m=0}^\infty \bigg[ \frac{\tilde r}{m+\tilde r_1}H\Bigr(0,0,\frac{\tilde r}{m+\tilde r_1}\Bigr)
-  \frac{\tilde r}{m+1}H\Bigr(0,0,\frac{-\tilde r}{m+1}\Bigr)\bigg].\label{CHdef}
\end{align}
\item[\rm(b)] The constant $C_H$ can be rewritten as follows:
\beq\ds
C_H=\sum_{k=1}^\infty b_{H,k}\hskip .3mm\mathcal{Q}_k(r_1,r_2),
\eeq
where $\mathcal{Q}_k(r_1,r_2)$'s are the constants defined in {\rm\eqnref{def:mu}} and $b_{H,k}$'s are the Taylor coefficients of $g(t):=H(t\mathbf{e}_3)-H(\mathbf{0})$, i.e., $H(t\mathbf{e}_3)=H(\mathbf{0})+b_{H,1} t +b_{H,2}t^2+\cdots,\quad t\in\RR$.
\end{itemize}
\end{theorem}
Note that we have $\tilde r=r/2$ and $\mu_1=\mu_2=\tilde{r}_1=\tilde{r}_2=1/2$ if $r_1=r_2=r$. Hence Eq.\;\eqnref{CHdef} coincides the series term in \eqnref{KLY1_a} for the case of two spheres of the same radius.

The second main result is the asymptotic formula for $\nabla u$, which shows the blow-up term explicitly in terms of elementary functions. We give the proof of Theorem \ref{main_thm2} in section \ref{subset:proof theorem}.
\begin{theorem}\label{main_thm2}
The solution $u$ to {\rm\eqnref{u:eqn}} admits the following decomposition in $ \mathbb{R}^3\setminus{\overline{(B_1\cup B_2)}}$:
$$
\nabla u(\Bx)= C_H \psi(\Bx) \left(\frac{\Bx-\Bp_1}{|\Bx-\Bp_1|^2}-\frac{\Bx-\Bp_2}{|\Bx-\Bp_2|^2}\right) + \nabla H(\Bx)+ r(\Bx),
$$
where \begin{gather}
\ds \psi(\Bx)=
\frac{\mu_\ep\tilde{r}}{2a_\ep}\bigg(\frac{\mu_1r_1}{|\Bx-\Bc_1|}+\frac{\mu_2\tilde{r}}{|\Bx-R_1(\Bc_2)|}
+\frac{\mu_2r_2}{|\Bx-\Bc_2|}+\frac{\mu_1\tilde{r}}{|\Bx-R_2(\Bc_1)|}\bigg),\end{gather}
and $\|r\|_\infty$ is bounded regardless of ${\ep}$.
\end{theorem}

\begin{cor}\label{cor:blowup}
Let $u$ be the solution to {\rm\eqnref{u:eqn}}. Then
\begin{center}$|\nabla u|$ blows up as $\ep$ tends to zero if and only if $\ds C_H\neq 0.$\end{center}
In particular, we have the followings:
\begin{itemize}
\item[\rm(a)] If $g=0$, then $|\nabla u|$ does not blow up.
\item[\rm(b)] If $g(t)=t^{2m-1}$ for some $m\in\NN$, then $|\nabla u|$ blows up.
\item[\rm(c)] If $r_1=r_2=r$ and $g$ is a polynomial of even degree, then $|\nabla u|$ does not blow up.
\end{itemize}
\end{cor}
\pf
Since $C_H$ is independent of $\ep$, Theorem \ref{main_thm2} asserts the equivalent condition for the gradient blow-up for $u$. Hence, we have (a). From \eqnref{polygamma}, we have $\psi_0(\tilde{r}_j)+\gamma<0$ for $j=1,2$. Hence $\mu_1$ and $\mu_2$ are positive.  We can similarly show that $\mathcal{Q}_{2m-1}>0$ for all $m\in\NN$, so that it follows (b).
If we assume $r_1=r_2=r$, then $\mathcal{Q}_{2m}(r,r)=0$. This proves (c).\qed


For example, if $H(\Bx)=b_2 p_2(\Bx)+b_3 p_3(\Bx)+b_4 p_4(\Bx)$ with $p_2(\Bx)=x_3^2-x_1^2$, $p_3(\Bx)=x_3^3-3x_3x_1^2,$ and $p_4(\Bx)=x_3^4-6x_1^2x_3^2+x_1^4$, then the corresponding electric field blows up if and only if $C_H=b_2 \mathcal{Q}_2+b_3\mathcal{Q}_3+b_4 \mathcal{Q}_4\neq 0$.
 Table \ref{table:nonlin} shows $\mathcal{Q}_k(r_1,r_2)$ for $r_1=1$ and various $r_2$ values.

\FloatBarrier
\begin{table}[!h]
\centering 
\renewcommand{\arraystretch}{1}
\setlength{\tabcolsep}{9pt}
\begin{tabular}{|c |r r r r r r |} 
\hline
$k$ & \multicolumn{1}{c}{$1$}& \multicolumn{1}{c}{$2$} & \multicolumn{1}{c}{$3$} &\multicolumn{1}{c}{ $4$}
&\multicolumn{1}{c}{$5$} &\multicolumn{1}{c|}{$6$}\\ [0.5ex] 
\hline \hline
$r_2=1.0$ &   $20.6709$ & \multicolumn{1}{c}{$0 $}      & $13.6009$ & \multicolumn{1}{c}{$0 $}     & $12.7843$ &\multicolumn{1}{c|}{$0$}\\
$r_2=0.7$ & $13.8369$& $-1.7996$  & 6.7967  &$-3.0177$ & 5.3858  & $-3.6830$ \\
$r_2=0.3$ &  3.9472  & $-1.1121$ &1.4317   &$-1.2751$ & 1.3062 & $-1.2938$ \\
$r_2=0.1$ & 0.5497& $-0.1795$ & 0.1851  & $-0.1828 $   & 0.1829 & $-0.1829$\\ 
\hline 
\end{tabular}\caption{values of the coefficients $\mathcal{Q}_k(r_1,r_2)$ when $r_1=1$} \label{table:nonlin}
\end{table}

\subsection{Superfocusing of the electric field}
\begin{figure}[!ht]
\begin{center}
\vskip -.1cm
\includegraphics[height=6cm]{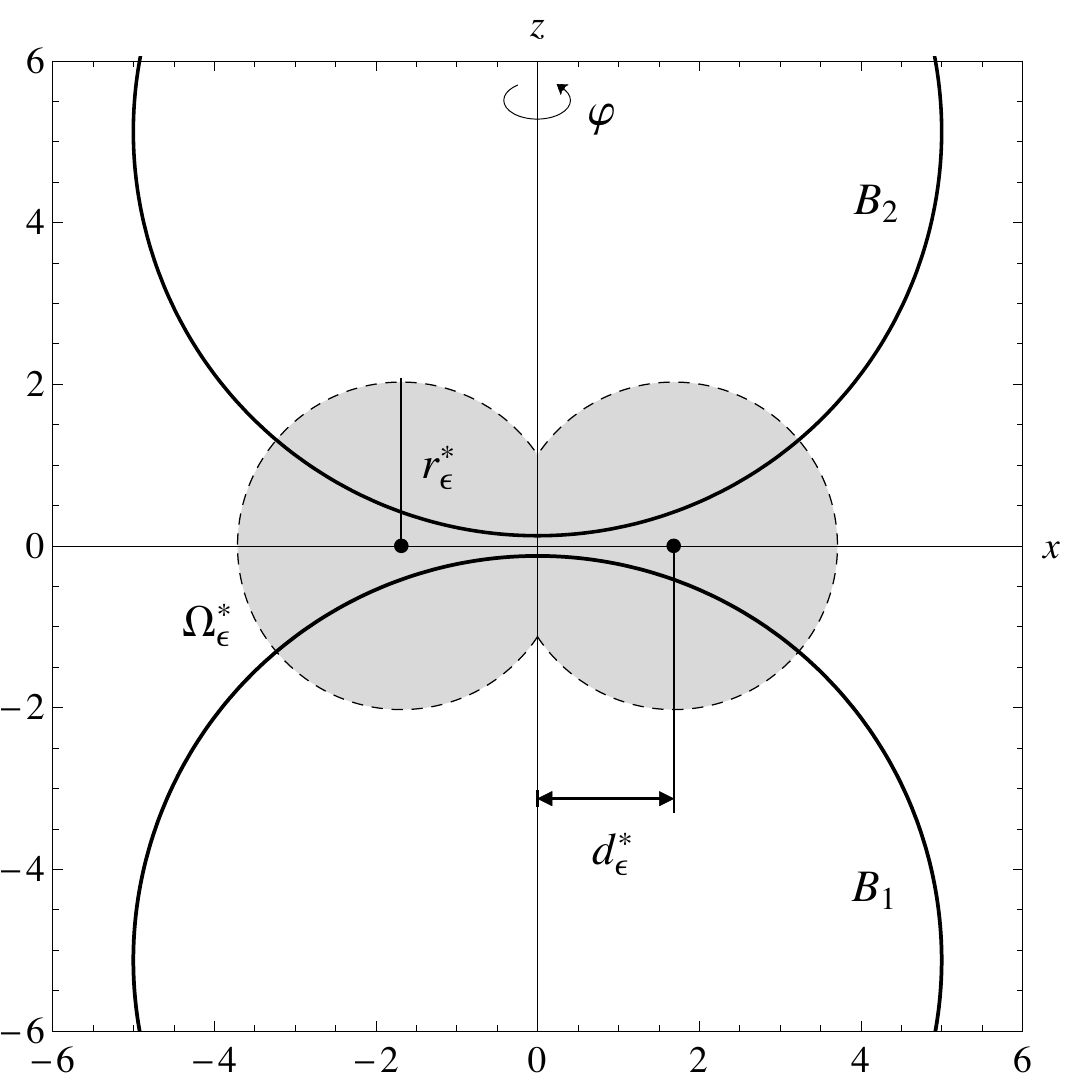}
\end{center}
\vskip -.5cm\caption{Superfocusing of the electric field: the gradient blow-up occurs only in the region $\Omega_\ep^*$ (shaded in the figure) whose dimension is of order $|\ln\ep|^{-1/2}$.}\label{fig:snowman}
\end{figure}

Let us denote $\theta_\ep=\sqrt{\ep|\ln\ep|}$ and
\beq\label{omegastar}
\Omega_\ep^*=\big\{\Bx\in\mathbb{R}^3:\big(|(x_1,x_2)|-d_\ep^*\big)^2+x_3^2 < \big(r_\ep^*\big)^2\big\}
\quad\mbox{with }
d_\ep^*=\frac{a_\ep}{\sin \theta_\ep},\ r_\ep^*={a_\ep}{\cot\theta_\ep},
\eeq
which is the rotation of the shaded region in Fig.\;\ref{fig:snowman} about the $x_3$-axis and can be written as $\Omega_{\ep}^*=\left\{\theta_\ep <\theta\leq\pi\right\}$ by the bispherical coordinates explained in the next section (see \eqnref{theta_const}). For small $\ep$, we have
$$
d_\ep^*,\;r_\ep^* \approx {\sqrt{2\tilde{r}}}\;{ |\ln\ep|^{-\frac{1}{2}}},$$
so that the width, length and height of $\Omega_\ep^*$ is of order $|\ln\ep|^{-1/2}$. This implies the convergence of the region $\Omega_\ep^*$ to the touching point $(0,0,0)$ as $\ep$ goes to $0$.

We show in the following theorem that the gradient blow-up occurs only in $\Omega_\ep^*$. In other words, we have superfocusing of the electric field confined in the narrow gap region between the two nearly touching metallic spheres. Moreover, the superfocusing region $\Omega_\ep^*$ is optimal in the sense that the order of its size cannot be smaller than $|\ln\ep|^{-1/2}$.
\begin{theorem}\label{prop:superfocusing}
To highlight the dependence on $\ep$, let us denote the solution to {\rm\eqnref{u:eqn}} by $u_\ep$. Then we have the followings.
\begin{itemize}
\item[\rm(a)]
The gradient blow-up occurs only in the region $\Omega_{\ep}^*$. More precisely, there exists a constant $C$ independent of $\ep$ satisfying
$$
|\nabla u_\ep (\Bx)| \leq C\quad \mbox{ for all } \Bx \in \mathbb{R}^3\setminus\overline{B_1\cup B_2\cup{\Omega_{\ep}^*}}.
$$
\item[\rm(b)]
The decaying order of $\theta_\ep$ is optimal in the following sense: if $C_H\neq0$ and $\tilde{\theta}_\ep$ satisfies $\lim_{\ep\rightarrow 0}\tilde{\theta}_\ep=0$ and $\lim_{\ep\rightarrow 0}\bigr({\tilde\theta_\ep}/{{\theta}_\ep}\bigr)=\infty,$ then we have
$$\inf_{x\in\tilde{\Omega}_\ep\setminus\overline{B_1\cup B_2}}|\nabla u_\ep(x)|\rightarrow \infty\quad\mbox{as }\ep\rightarrow 0,$$ where $\tilde{\Omega}_\ep$ is defined as {\rm\eqnref{omegastar}} with $\tilde{\theta}_\ep$ in the place of $\theta_\ep$.
\end{itemize}
\end{theorem}
We will prove the theorem in section \ref{subset:proof theorem}.

\section{Bispherical coordinate system}\label{section:spherical}
\begin{figure}[!ht]
\begin{center}
\vskip -.1cm
\includegraphics[height=7cm]{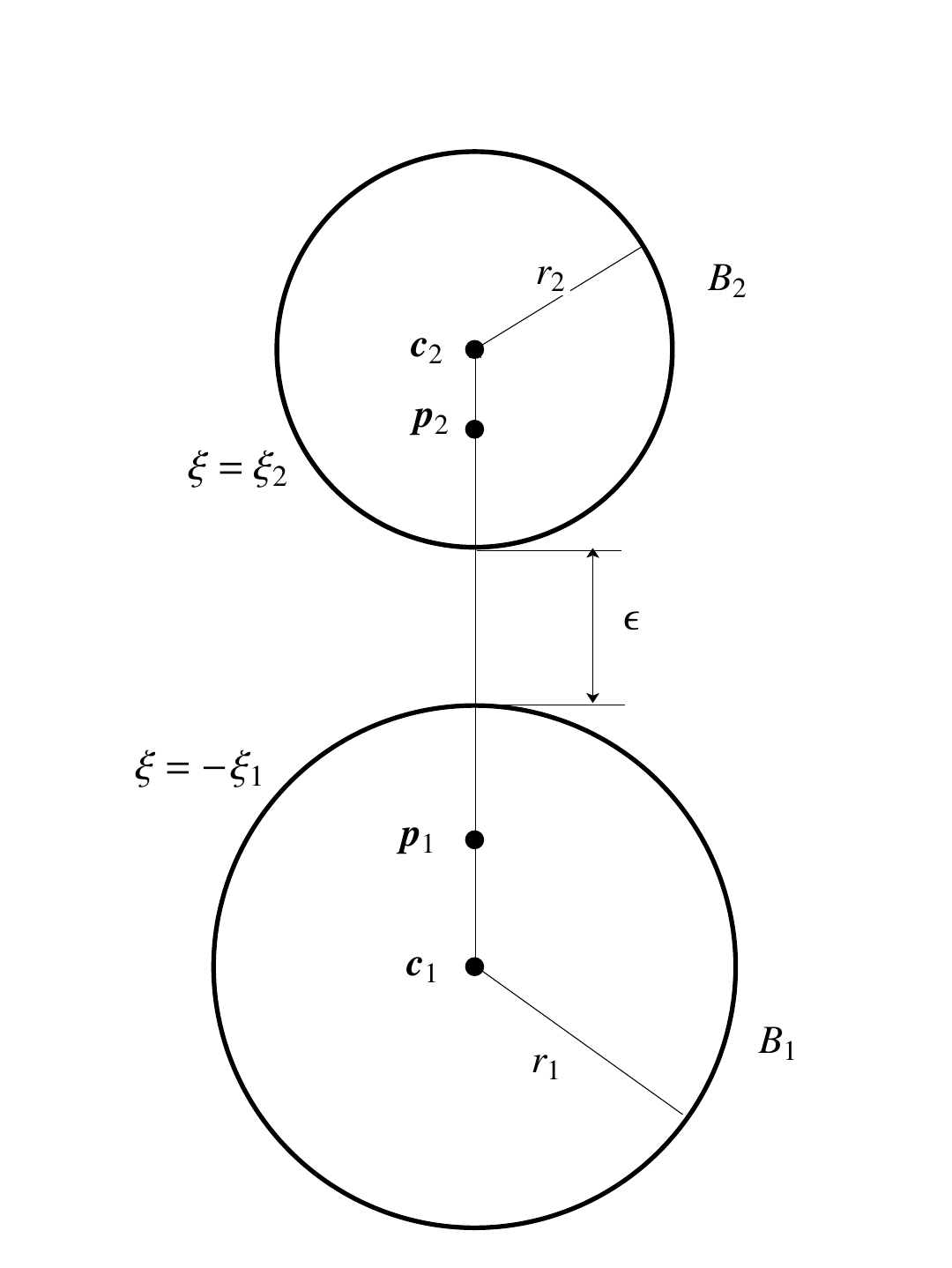}\hskip 7mm
\includegraphics[height=6.5cm]{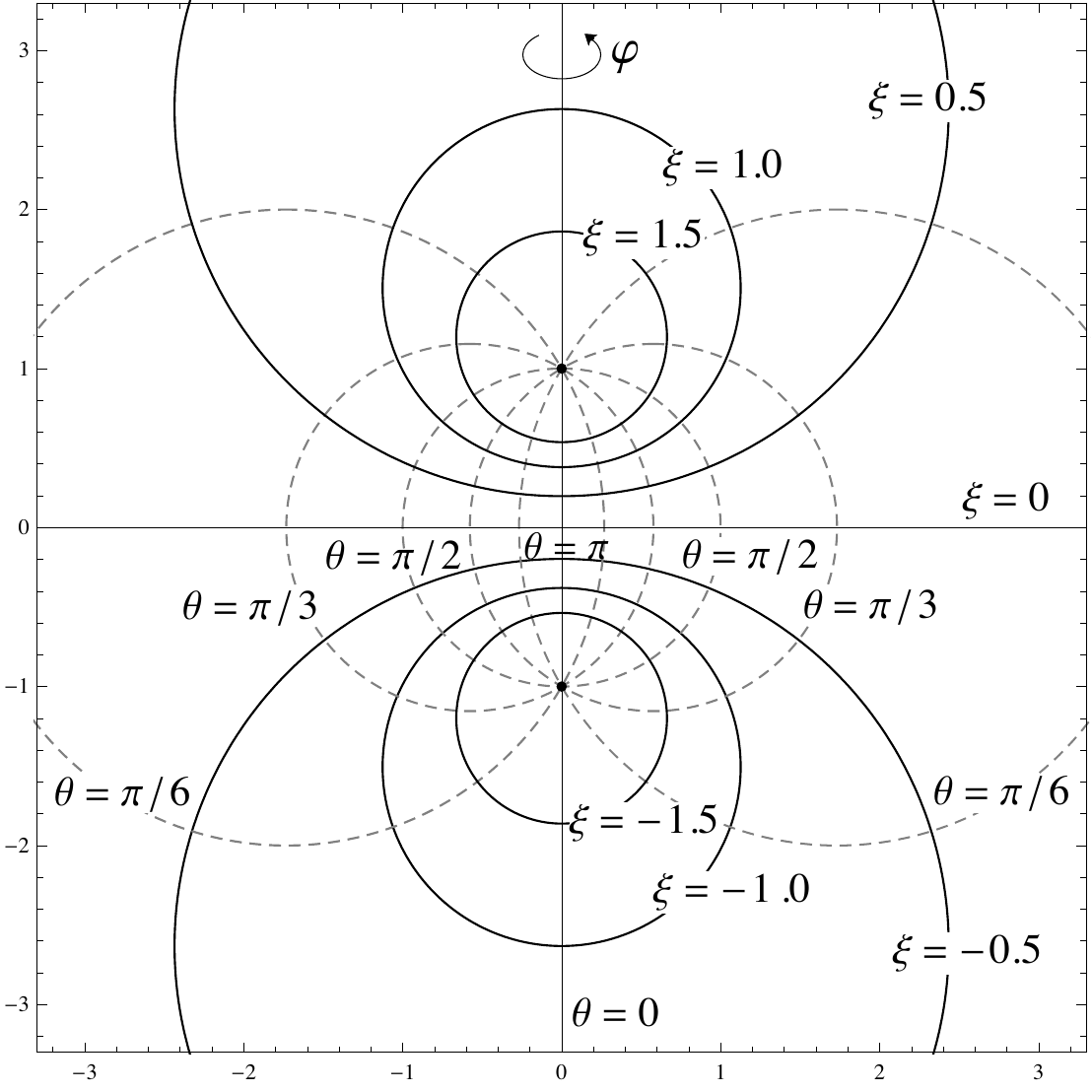}
\end{center}
\vskip -.5cm\caption{A pair of spherical perfect conductors (the left) and the coordinate level curves for the bipherical coordinate system with $a_\ep=1$ (the right).}\label{fig:fig_nonuniform}
\end{figure}

Let us introduce the bispherical coordinate system $(\xi,\theta,\varphi)\in\RR\times[0,\pi]\times[0,2\pi)$ with poles located at $\Bp_1$ and $\Bp_2$.
Each $\Bx=(x_1,x_2,x_3)$ in the Cartesian coordinate system of $\RR^3$ corresponds to $(\xi,\theta,\varphi)$ through
\beq\label{bipolar}e^{\xi-i\theta}=\frac{z+a_\ep}{z-a_\ep}\quad\mbox{with }z=x_3+i\left|(x_1,x_2)\right|\eeq
with $\varphi$ the angle of rotation about the $x_3$-axis. One can rewrite the Cartesian coordinates in terms of the bispherical coordinates as \begin{align*}
\ds x_1=a_\ep\frac{\sin\theta\cos\varphi}{\cosh\xi-\cos\theta},\quad
\ds x_2=a_\ep\frac{\sin\theta\sin\varphi}{\cosh\xi-\cos\theta},\quad
\ds x_3=a_\ep\frac{\sinh\xi}{\cosh\xi-\cos\theta}.
\end{align*}
It can be easily shown that the coordinate surfaces $\{\xi=c\}$ and $\{\theta=c\}$ for a nonzero $c$ are respectively the zero level set of
\begin{align}\ds f^\xi(x_1,x_2,x_3)&=\left(x_3-a_\ep\coth c\right)^2 +|(x_1,x_2)|^2-\left(\frac{a_\ep}{\sinh c}\right)^2,\label{fxi}\\
\label{theta_const}
\ds f^\theta(x_1,x_2,x_3)&=\big(\left|(x_1,x_2)\right|-{a_\ep}{\cot c}\big)^2+x_3^2 -\left(\frac{a_\ep}{\sin c}\right)^2.\end{align}
We illustrate the coordinate surfaces of the bispherical coordinate in Fig.\;\ref{fig:fig_nonuniform}.

Note that $\Bp_1$ and $\Bp_2$ are contained in $\{x_2=0\}$ and they are again fixed points of the combined reflections w.r.t. the two circles $\p B_j\cap\{x_2=0\}$, $j=1,2$.
Remind the circle of Apollonius in two dimensions: for disk $B_r(\Bc)$ and $\mathbf{y}\notin\overline{B_r(\Bc)}$, the circle $\p B_r(\Bc)$ is the locus of $\Bx$ satisfying
\beq\notag
\frac{|\Bx-\mathbf{y}|}{|\Bx-R(\mathbf{y})|}=\frac{|\mathbf{y}-\Bc|}{r},\eeq
where $R$ is the reflection w.r.t. the circle $\p B_r(\Bc)$. Applying this property to $\p B_j\cap\{x_2=0\}$, we have that ${|\Bx-\Bp_1|}/{|\Bx-\Bp_2|}$ is constant on $\p B_j\cap\{x_2=0\},\ j=1,2.$ From \eqnref{bipolar} and the fact
$$\left|\frac{x_3+i|x_1|+a_\ep}{x_3+i|x_1|-a_\ep}\right|=\frac{|\Bx-\Bp_1|}{|\Bx-\Bp_2|},$$
$\xi$ is constant on $\p B_j\cap\{x_2=0\}$, and hence so does on $\p B_j$. Furthermore, applying Eq.\;\eqnref{fxi}, we obtain for $j=1,2$ that
\begin{gather}\label{partialBj}
\ds\p B_j = \big\{\xi=(-1)^j\xi_j\big\}\quad\mbox{and}\quad
\Bc_j=(-1)^{j}a_\ep\coth\xi_j\mathbf{e}_3
\end{gather}
with two positive constants $\xi_1$ and $\xi_2$ given by
\beq\label{xi0909}\xi_j=\sinh^{-1}\left({a_\ep}{r_j^{-1}}\right)=\alpha_\ep r_j^{-1}+O(\ep\sqrt{\ep}).\eeq
Eq.\;\eqnref{def:a} and Eq.\;\eqnref{xi0909} imply
\beq\label{axi1xi2} \frac{a_\ep}{\xi_1+\xi_2}=\tilde{r}+O({\epsilon})\eeq
and
 \beq\label{xi3range}
\begin{cases}-\xi_1\leq\xi\leq\xi_2\quad &\mbox{in }\RR^3\setminus\overline{(B_1\cup B_2)}\\
\xi\in(-\infty, -\xi_1)\quad&\mbox{in }B_1\\
\xi\in (\xi_2,\infty)\quad&\mbox{in } B_2.
\end{cases}\eeq

Let us consider the multiple reflections of the spheres centers $\Bc_1$ and $\Bc_2$ w.r.t. the two spheres.
From the definition of the reflection $R_1$ and \eqnref{partialBj}, it follows
\begin{align*}
\mathbf{e}_3\cdot R_1(\Bc_2)&=-a_\ep\coth\xi_1+\frac{r_1^2}{|\Bc_1-\Bc_2|}\\
&=-a_\ep\frac{\cosh\xi_1}{\sinh\xi_1}+a_\ep\frac{1}{\sinh^2\xi_1}\frac{1}{\frac{\cosh\xi_1}{\sinh\xi_1}+\frac{\cosh\xi_2}{\sinh\xi_2}}\\
&=-a_\ep\coth(\xi_1+\xi_2).
\end{align*}
By the same way for $m=0,1,\dots$, we have
\beq\label{refcenter}
\begin{cases}
\ds(R_1\circ R_2)^m(\Bc_1)=-\mathbf{p}_m^{\xi_1},\quad
\ds(R_2\circ R_1)^m(\Bc_2)=\mathbf{p}_m^{\xi_2}, \\[2mm]
\ds(R_2\circ R_1)^m\circ R_2(\Bc_1)=-(R_1\circ R_2)^m\circ R_1(\Bc_2)=\mathbf{p}_m^{\xi_1+\xi_2},
\end{cases}\eeq
where
\beq\label{eqnpm}\mathbf{p}_m^c=a_\ep\coth\bigr(m(\xi_1+\xi_2)+c\bigr)\mathbf{e}_3\quad\mbox{for }c=\xi_1,\xi_2,\xi_1+\xi_2.\eeq

\subsection{Scale factors and harmonic functions}
The bispherical coordinate system $(\xi,\theta,\varphi)$ is an orthogonal coordinate system. We denote its orthogonal coordinate directions as $\{\hat{\mathbf{e}}_{\xi},\hat{\mathbf{e}}_{\theta}, \hat{\mathbf{e}}_{\varphi}\}$, {\it i.e.},
\beq\label{basis}
\hat{\mathbf{e}}_{\xi}= \frac{\p \mathbf{x}/ \p \xi}{|\p \mathbf{x}/ \p\xi|}, \quad  \quad \hat{\mathbf{e}}_{\theta}= \frac{\p \mathbf{x}/ \p \theta}{|\p \mathbf{x}/ \p\theta|},
\quad \quad \hat{\mathbf{e}}_{\varphi}= \frac{\p \mathbf{x}/ \p \varphi}{|\p \mathbf{x}/ \p\varphi|}.
\eeq
The scale factors for the bispherical coordinates are
\beq \label{scalef}\sigma_\xi=\sigma_\theta=\frac{a_\ep}{\cosh \xi-\cos\theta}\quad\mbox{and}\quad \sigma_\varphi = \frac{a_\ep\sin\theta}{\cosh \xi-\cos\theta},\eeq
so that the gradient for scalar valued function $g$ can be written as
\beq \label{grad_bipolar}
\nabla g = \frac{1}{\sigma_\xi}\frac{\p g}{\p\xi}\hat{\mathbf{e}}_{\xi}+ \frac{1}{\sigma_\theta} \frac{\p g}{\p\theta}\hat{\mathbf{e}}_{\theta}+\frac{1}{\sigma_\varphi }\frac{\p g}{\p\varphi}\hat{\mathbf{e}}_{\varphi}.
\eeq
Here and in the remaining of the paper, the symbol $\nabla$ denotes the gradient in the Cartesian coordinates.
It can be also shown
\beq\label{exi:x}
\hat{\mathbf{e}}_{\xi}(\Bx) = \sigma_\xi(\Bx)\mathbf{N}(\Bx)\quad\mbox{with }\mathbf{N}(\Bx)=\left(\frac{\Bx-\Bp_1}{|\Bx-\Bp_1|^2}-\frac{\Bx-\Bp_2}{|\Bx-\Bp_2|^2}\right).\eeq

As one can see in Page 111 of \cite{MS}, any harmonic function $f$ has a general  $R$-separation
\begin{align}\ds f(\xi,\theta,\varphi)&=\sqrt{\cosh\xi-\cos\theta}\sum_{n=1}^{+\infty}\sum_{m=0}^n \Big[D_n^m e^{(n+\frac{1}{2})|\xi|}+E_n^m e^{-(n+\frac{1}{2})|\xi|}\Big]\notag\\\label{eqn:SVsolution}
\ds&\qquad\times P_n^m(\cos\theta)\Big[F_n^m \cos(m\varphi)+G_n^m\sin(m\varphi)\Big],
\end{align}
where $P_n^m$'s are the Legendre associated functions and $D_n^m$, $E_n^m,\ F_n^m$, and $G_n^m$ are constants.
It is well known that the generating function for the Legendre polynomials $P_n(x)$'s, which are $P_n^0(x)$'s, is given by
\beq\notag
\frac{1}{\sqrt{1-2xt+t^2}}=\sum_{n=0}^\infty t^n P_n(x) \quad\mbox{for }|x|\leq 1, |t|< 1,
\eeq
and they form an orthogonal basis of $L^2[0,1]$. From the equation above, the constant function $1$ can be expressed as
\begin{align}\label{ident_const}
\ds 1=\sqrt{2}\sqrt{\cosh\xi-\cos\theta} &\sum_{n=0}^\infty e^{-\left(n+\frac{1}{2}\right)|\xi|}P_n(\cos\theta).\end{align}
We also have the following two identities for $\xi\in\RR$ \cite{MS}: \begin{align}
\label{Legen_int1}
\ds\int_{-1}^1 \frac{P_n(s)}{(\cosh\xi-s)^{\frac 1 2} }~ ds&=\frac{2\sqrt{2}}{2n+1}e^{-(n+\frac{1}{2})|\xi|}
,\\[1mm]\label{Legen_int2}
\ds\int_{-1}^1 \frac{P_n(s) }{(\cosh\xi-s)^{\frac 3 2}}~ds&=\frac{2\sqrt{2}}{\sinh|\xi|}e^{-(n+\frac{1}{2})|\xi|}.
\end{align}


\section{The singular function $h$}
In this section, we give the two series expansions for the solution $h$ to \eqnref{h:eqn} by the bispherical coordinates and give the proof of Theorem \ref{lem:CH}.

\subsection{Solution by separation of variables}\label{sol_sepa}
We set \begin{align}
\ds C_j=\frac{(-1)^j}{8\pi a_\ep}  \frac{U(\xi_j)-U(0)}{U(\xi_1)U(\xi_2)-U^2(0)},\quad j=1,2,
\label{def_C}
\end{align}
where
 \begin{align} U(c)=\sum_{n=0}^\infty\frac{e^{(2n+1)c}}{e^{(2n+1)(\xi_1+\xi_2)}-1}\quad\mbox{for }  0<c<\xi_1+\xi_2.\label{def:U}\end{align}
In the following lemma we express $h$ in the form of \eqnref{eqn:SVsolution} with the coefficients defined using $C_j$'s, where $C_j$\rq{}s are actually potential values of $h$ on $\p B_j$'s. Note that $h$ is independent of $\varphi$ due to its symmetry under the rotation about $x_3$-axis. We omit the variable $\varphi$ in $h$ for notational simplicity.
\begin{lemma}\label{h_series_expand}
The solution $h$ to {\rm \eqnref{h:eqn}} can be represented as
\begin{align}
\ds h(\xi,\theta)=\sqrt{2}\sqrt{\cosh\xi-\cos\theta}\sum_{n=0}^\infty \left( A_n e^{(n+\frac{1}{2})\xi}
+
B_ne^{-(n+\frac{1}{2})\xi} \right)P_n(\cos\theta),\quad -\xi_1\leq\xi\leq\xi_2, \label{h_series}
\end{align}
where
$$
A_n=\frac{C_2e^{(2n+1)\xi_1}-C_1}{e^{(2n+1)(\xi_1+\xi_2)}-1} \quad\mbox{and}\quad
B_n=\frac{C_1e^{(2n+1)\xi_2}-C_2}{e^{(2n+1)(\xi_1+\xi_2)}-1}.
$$
Moreover, $h$ satisfies\beq\label{eqn:Lekner}
h|_{\p B_j}=C_j, \quad j=1,2.
\eeq
\end{lemma}
\pf
Let us denote the right-hand side of \eqnref{h_series} as $\tilde{h}$. In the following, we prove $\tilde{h}$ satisfies all the constraints in \eqnref{h:eqn}.

One can easily show that all the terms in the series expansion of $\tilde h$ are harmonic, see \eqnref{eqn:SVsolution}, and they are exponentially decay (uniformly for $\xi\in[-\xi_1,\xi_2]$) in $n$. Hence, $\tilde h$ is harmonic. From \eqnref{ident_const}, $\tilde h$ is constant on $\p B_1$ and $\p B_2$. More precisely,
$$
\tilde{h}(-\xi_1,\theta)=C_1\quad\mbox{and}\quad \tilde{h}(\xi_2,\theta)=C_2.
$$

It can be easily shown that the outward unit normal vector $\nu$ to $\p B_j\left(=\{\xi=(-1)^j\xi_j\}\right)$ is
\beq\label{eqn:nuexi}\nu={(-1)^{j+1}}\hat{\mathbf{e}}_{\xi}
\eeq
and a sufficiently smooth  function $v$ satisfies
\beq\label{surf_int}
\int_{\p B_j} \partial_\nu v ~d\sigma = (-1)^{j+1}\int_0^{2\pi}\int_0^\pi \left(  \frac{\p v }{\p\xi}\bigg|_{\xi=(-1)^j\xi_j} \sigma_\varphi(\xi_j,\theta,\varphi) \right) d\theta d\varphi.
\eeq
In particular, using \eqnref{Legen_int1}, we derive for $j=1,2$ that
\begin{align*}
\ds
\sqrt{2}\int_{\p B_j}{\p_\nu}\left(\sqrt{\cosh\xi -\cos\theta}\ e^{(n+\frac{1}{2})\xi}P_n(\cos\theta)\right) d\sigma&=-8\pi a_\ep\delta_{2j},\\
\ds
\sqrt{2}\int_{\p B_j}\p_\nu\left(\sqrt{\cosh\xi -\cos\theta}\ e^{-(n+\frac{1}{2})\xi}P_n(\cos\theta)\right) d\sigma&=-8\pi a_\ep\delta_{1j},
\end{align*}
where $\delta_{i,j}$ is 1 if $i=j$ and zero otherwise.
Hence we have
$$
\int_{\p B_j}\p_\nu \tilde{h} ~d\sigma = -8\pi a_\ep\left(\sum A_n\right)\delta_{2j}-8\pi a_\ep\left(\sum B_n\right)\delta_{1j}=(-1)^{j+1} \quad \mbox{for }j=1,2.
$$

Now it only remains to show the decay property at infinity for $\tilde{h}$. In fact, it is enough to show that $\tilde{h}(\Bx)=O(|\Bx|^{-1})$ as $|\Bx|\rightarrow\infty$ because the total flux on $\p B_1\cup \p B_2$ is zero.
Note that the radial distance $|\Bx|$ satisfies
\beq\label{rinf}|\Bx|=a_\ep\sqrt{\frac{\cosh\xi+\cos\theta}{\cosh\xi-\cos\theta}},\eeq
so that $|\Bx|\to\infty$ if and only if $(\xi,\theta)\to(0,0)$.
Hence we only need to show
\beq\notag
\limsup_{(\xi,\theta)\rightarrow(0,0)} \frac{|\tilde{h}(\xi,\theta)|}{\sqrt{\cosh\xi-\cos\theta}}\leq C,
\eeq
for some constant $C$ independent of $\xi$ and $\theta$.
Owing to
$$
\frac{|\tilde{h}(\xi,\theta)|}{\sqrt{\cosh\xi-\cos\theta}}
\leq \sum_{m=0}^\infty \left(|A_n|e^{(n+\frac{1}{2})\xi_2} +|B_n|e^{(n+\frac{1}{2})\xi_1}\right) <\infty\quad\mbox{for all }-\xi_1\leq\xi\leq\xi_2,
$$
the decay condition follows. This completes the proof.
\qed

The following asymptotic of $U$ has been derived by J. Lekner in \cite{lekner11}.
\begin{lemma}{\rm(\cite{lekner11})}
For small $\ep>0$, the function $U$ defined in {\rm\eqnref{def:U}} satisfies
$$U(\xi_j)=\frac{1}{2(\xi_1+\xi_2)}\left[\ln\Big(\frac{2}{\xi_1+\xi_2}\Big)
-\psi_0\Big(1-\frac{\xi_j}{\xi_1+\xi_2}\Big)\right]+O(\sqrt{\ep}),$$
and
$$
U(0)=\frac{1}{2(\xi_1+\xi_2)}\left[\ln\Big(\frac{2}{\xi_1+\xi_2}\Big)
+\gamma\right]+O(\sqrt{\ep}).
$$
\end{lemma}

\begin{cor}\label{lemma:lekner}
We have\beq C_j=(-1)^j\mu_\ep\hskip .3mm\mu_j+O(\ep),\quad j=1,2,\label{h_const} \eeq
with $\mu_\ep,\mu_j$'s given in {\rm\eqnref{mu1002}}.
\end{cor}


\subsection{Expansion by Potentials of point charges}\label{sol_lpm}
The fundamental solution $\Gamma$ to the Laplacian in three dimensions is given by
$$\Gamma(\Bx)=-\frac{1}{4\pi|\Bx|}.$$
We can rewrite $\Gamma$ by the bispherical coordinates as a fraction of $w_\theta$ which is defined as
 \beq\label{def:w_theta} w_\theta(\xi):=\sqrt{\cosh\xi-\cos\theta}.\eeq
\begin{lemma}\label{def_g2}
Let $\Bx=(x_1,x_2, x_3)\in\RR^3$ be a point with the bispherical coordinate $(\xi,\theta,\varphi)$.
Then for $\xi\in\RR$ we have
$$
\frac{w_\theta(\xi)}{w_\theta(\xi-2\xi_0)}
=-\frac{4\pi a_\ep}{\sinh|\xi_0|}\Gamma(\Bx-\Bx_0)\quad \mbox{with }\Bx_{0}= a_\ep\coth \xi_0 \mathbf{e}_3.$$

\end{lemma}
\pf
We have from \eqnref{bipolar} that
$$ x_3+i|(x_2,x_3)|=z=\frac{2a_\ep}{e^{\xi-i\theta}-1}+a_\ep.$$
Note that $$\coth\xi_0=\frac{\sinh 2\xi_0}{\cosh 2\xi_0 -1}=\frac{2}{e^{2\xi_0}-1}+1\quad\mbox{and}\quad|e^{\xi-i\theta}-1|=\sqrt{2e^\xi (\cosh\xi-\cos\theta)}. $$
 Hence it follows
\begin{align*}
|\Bx-\Bx_0|&=\Big|(x_1,x_2,x_3-a_\ep\coth \xi_0)\Big|=\Big|x_3+i|(x_2,x_2)|-a_\ep\coth(\xi_0)\Big|\\
&=2a_\ep\left|\frac{1}{e^{\xi-i\theta}-1}-\frac{1}{e^{2\xi_0}-1}\right|
=2a_\ep\left|\frac{e^{2\xi_0}(e^{\xi-2\xi_0-i\theta}-1)}{(e^{2\xi_0}-1)(e^{\xi-i\theta}-1)}\right|=\frac{a_\ep}{\sinh|\xi_0|}\frac{w_\theta(\xi-2\xi_0)}{w_\theta(\xi)}.\end{align*}
Thus we prove the lemma.
\qed


\smallskip
Let us denote $S$ for $(\xi,\theta,s)\in\RR\times[0,\pi]\times\RR$ as
\begin{align}\label{def_S}
\ds S(\xi,\theta;s):&=\sum_{m=0}^\infty\frac{w_\theta(\xi)}{w_\theta\bigr(\xi-2m(\xi_1+\xi_2)-s\bigr)}.
\end{align}
Because $w_\theta$ is an even function for $\xi$ and $$\frac{w_\theta(\xi)}{w_\theta(\xi\pm 2\xi_0)}
=-\frac{4\pi a_\ep}{\sinh|\xi_0|}\Gamma(\Bx\pm\Bx_0), \quad\Bx_{0}= a_\ep\coth \xi_0 \mathbf{e}_3,$$
it follows
\begin{align}\label{SGamma}
 \ds S(\pm\xi,\theta;2c)&=-\sum_{m=0}^\infty \frac{4\pi a_\ep}{\sinh\left|\xi_m^{c}\right|}{\Gamma\left(\Bx\mp a_\ep \coth(\xi_m^{c})\mathbf{e}_3\right)},\quad \xi_m^{c}=m(\xi_1+\xi_2)+c.\end{align}
We can express the solution $h$ to \eqnref{h:eqn} as a linear combination of $S$.
\begin{lemma}\label{h_L}
We have
\begin{align*}
\ds h(\Bx)&=C_2S\bigr(\xi,\theta;2\xi_2\bigr)-C_1S\big(\xi,\theta;2(\xi_1+\xi_2)\big)\\
 \ds &+C_1S\big(-\xi,\theta;2\xi_1\big)-C_2S\big(-\xi,\theta;2(\xi_1+\xi_2)\big).
\end{align*}
\end{lemma}
\pf
Since $\xi_1,\xi_2>0$, we have
$$
\frac{1}{e^{(2n+1)(\xi_1+\xi_2)}-1}=\sum_{m=0}^\infty e^{-(m+1)(2n+1)(\xi_1+\xi_2)}
=\sum_{m=0}^\infty e^{-(n+\frac{1}{2})[2m(\xi_1+\xi_2)+2(\xi_1+\xi_2)]}.
$$
Applying the above identity and interchanging the order of summation which is possible due to the absolute convergence of the series, Eq.\;\eqnref{h_series} becomes
\begin{align}
\ds{h(\Bx)}\notag
\ds&={\sqrt{2}\sqrt{\cosh\xi-\cos\theta}}\\ \ds&\times\sum_{m=0}^\infty
\sum_{n=0}^\infty \bigg[ C_2 e^{-(n+\frac{1}{2})\xi_{m,1}}-C_1e^{-(n+\frac{1}{2})\xi_{m,2}}
+C_1e^{-(n+\frac{1}{2})\xi_{m,3}}-C_2e^{-(n+\frac{1}{2})\xi_{m,4}}\bigg]P_n(\cos\theta)\notag
\end{align}
with \begin{align*}
\ds \xi_{m,1}&= -\xi+{2m(\xi_1+\xi_2)+2\xi_2},\quad \xi_{m,2}=-\xi+{2m(\xi_1+\xi_2)+2(\xi_1+\xi_2)},\\
\ds\xi_{m,3}&=\xi+{2m(\xi_1+\xi_2)+2\xi_1}, \quad \xi_{m,4}=\xi+{2m(\xi_1+\xi_2)+2(\xi_1+\xi_2)}.
\end{align*}
Thanks to \eqnref{ident_const}, we obtain
\begin{align*}\sum_{n=0}^\infty e^{-(n+\frac{1}{2})\xi_{m,j}}P_n(\cos\theta)&=\frac{1}{\sqrt{2}}\frac{1}{w_\theta(\xi_{m,j})},\quad j=1,\dots,4.
\end{align*}
This completes the proof.
\qed

\begin{cor}\label{h_Gamma}
For $c$ is either $\xi_1,\xi_2$ or $\xi_1+\xi_2$, we denote
\beq\label{pmqm}
\xi_m^c=m(\xi_1+\xi_2)+c,\quad q_m^c=\frac{4\pi a_\ep}{ \sinh\xi_m^c},\quad \mathbf{p}_m^c=a_\ep\coth\xi_m^c\mathbf{e}_3.
\eeq
Then we can expand $h$ and $C_H^\ep$ as
\begin{align}\label{hasymp}
\ds h(\Bx)&=
C_1\sum_{m=0}^\infty \Bigr[ q_m^{\xi_1+\xi_2} \Gamma(\Bx-\mathbf{p}_m^{\xi_1+\xi_2})
-q_m^{\xi_1} \Gamma(\Bx+\mathbf{p}_m^{\xi_1})\Bigr]\notag\\
& \quad -
\ds C_2\sum_{m=0}^\infty \Bigr[ q_m^{\xi_2} \Gamma(\Bx-\mathbf{p}_m^{\xi_2})
-q_m^{\xi_1+\xi_2} \Gamma(\Bx+\mathbf{p}_m^{\xi_1+\xi_2})\Bigr]
\end{align}
and
\begin{align}
\ds C_H^{\ep} &=\frac{C_1}{C_1-C_2}\sum_{m=0}^\infty \Bigr[ q_m^{\xi_1+\xi_2} H\left(\mathbf{p}_m^{\xi_1+\xi_2}\right)
-q_m^{\xi_1} H\left(-\mathbf{p}_m^{\xi_1}\right)\Bigr]\notag\\
&\quad
\ds-\frac{C_2}{C_1-C_2}\sum_{m=0}^\infty \Bigr[ q_m^{\xi_2} H\left(\mathbf{p}_m^{\xi_2}\right)
-q_m^{\xi_1+\xi_2} H\left(-\mathbf{p}_m^{\xi_1+\xi_2}\right)\Bigr].\label{CHexp}
\end{align}
\end{cor}
\pf
From \eqnref{SGamma} and Lemma \ref{h_L}, we prove \eqnref{hasymp}. Remind that \eqnref{eqn:Lekner} implies
$$h|_{\p B_1}-h|_{\p B_2}=C_1-C_2.$$ One can easily show \eqnref{CHexp} by computing $(u|_{\p B_1}-u|_{\p B_2})$ from \eqnref{0417:udiff} and \eqnref{hasymp}.
\qed

Note that $\mathbf{p}_m^c$'s are multiply reflected points of $\Bc_1$ and $\Bc_2$ with respect to the two spheres, see \eqnref{refcenter}.
Corollary \ref{h_Gamma} has the same formality as Lemma 4.1 in \cite{LY}, where a recursively defined series was used in the place of $q_m^c$.
In this paper, we are able to have formulas much simpler than those in \cite{LY} thanks to adopting the bispherical coordinate system.

\subsection{Proof of Theorem \ref{lem:CH}} \label{section:CH}
 \begin{lemma}\label{CHC}
Let $\mathbf{p}_m^c,q_m^c,\xi_m^c$ be given as in Corollary {\rm \ref{h_Gamma}} and $H$ be an entire harmonic function.
Then there is a constant $C$ independent of $\ep$ satisfying
\beq\notag
\bigg|\sum_{m=0}^\infty q_m^c H(\mathbf{p}_m^c)-\sum_{m=0}^\infty \frac{4\pi a_\ep}{\xi_m^c} H\Big(0,0,\frac{a_\ep}{\xi_m^c}\Big) \bigg| \leq{C}\ep|\ln\ep|.
\eeq

\end{lemma}
\pf
Let $c$ be fixed to be either $\xi_1$, $\xi_2$ or $\xi_1+\xi_2$.
We can assume $H(\mathbf{0})=0$ since the constant term of $H$ does not change the gradient of the potential function. Then $H$ can be written as
$$H(\Bx)=\sum_{k=1}^\infty H_k(\Bx),$$
where $H_k$'s are homogeneous polynomials in $\Bx$ of degree $k$. Especially, we have
$$H(x_k\mathbf{e}_3)=\sum_{k=1}^\infty H_k(x_3\mathbf{e}_3)=\sum_{k=1}^\infty b_{H,k}x_3^k\quad\mbox{with }b_{H,k}= \frac{1}{k!}\frac{\p^k H}{\p x_3^k}(\mathbf{0}).$$

From \eqnref{pmqm}, one obtains
$$q_m^c H_k(\mathbf{p}_{m}^c)
=
{4\pi b_{H,k} a_\ep^{k+1} }
\frac{\cosh^k\big(s_0m+c\big) }{ \sinh^{k+1}\big(s_0 m+c\big)},\quad s_0=\xi_1+\xi_2.$$
We denote
\begin{align}\notag
G_k&=4\pi b_{H,k}a_\ep^{k+1} \sum_{m=0}^\infty \frac{1}{(s_0m+c)^{k+1}},\\
\label{fk}f_k(x)&=\frac{\cosh^k x}{\sinh^{k+1} x}-\frac{1}{x^{k+1}}, \qquad x>0,~k\geq1 .
\end{align}
Then it follows
\begin{align}\notag
 R_k:=\sum_{m=0}^\infty q_m^c H_k(\mathbf{p}_m^c)-G_k
&= \frac{b_{H,k} 4\pi a_\ep^{k+1}}{s_0} \sum_{m=0}^\infty f_k(s_0 m+c) s_0.
\end{align}

Thanks to Lemma \ref{Sum_fk_lem} in section \ref{section:euler}, one obtains that there is a constant $C$ independent of $\ep$ and $k$ such that
\beq\notag
\left| R_k\right| \leq C \ep|\ln\ep||b_{H,k}| 2^k \quad \mbox{for all } k\geq1.\eeq
Note that the series $\sum_{k=1}^\infty |b_{H,k}| 2^k$ converges since $H$ is an entire function, so we have
 \begin{align}\label{ineq:0717}\bigg|\sum_{m=1}^\infty \sum_{k=1}^\infty q_m^c H_k(\mathbf{p}_m^c)-\sum_{k=1}^\infty G_k\bigg|&=\bigg|\sum_{k=1}^\infty R_k\bigg|\leq \tilde{C}\ep|\ln\ep|
 \end{align}
for a constant $\tilde C$ independent of $\ep$ and $k$. The first equality in \eqnref{ineq:0717} holds because of the absolute convergence of two series in the leftmost side. Because of the same reason, we also have
\begin{gather*}\sum_{m=1}^\infty \sum_{k=1}^\infty q_m^c H_k(\mathbf{p}_m^c)=\sum_{m=0}^\infty q_m^c H(\mathbf{p}_m^c),\\
\sum_{k=1}^\infty G_k =4\pi \sum_{m=0}^\infty\sum_{k=1}^\infty b_{H,k} \frac{a_\ep^{k+1}}{(s_0m+c)^{k+1}}=\sum_{m=0}^\infty \frac{4\pi a_\ep}{\xi_m^c} H\Big(0,0,\frac{a_\ep}{\xi_m^c}\Big) .
\end{gather*}
Therefore we prove the theorem thanks to \eqnref{ineq:0717}.
\qed

{{\noindent{{\textbf{Proof of Theorem \ref{lem:CH}}}}}
Thanks to \eqnref{def:a} and \eqnref{def:tilder}, one can easily show
\beq\notag
\quad \frac{a_\ep}{\xi_m^c}=
\begin{cases}
\ds \frac{\tilde{r}}{m+1}+O\Bigr(\frac{\ep}{m+1}\Bigr) \quad &\mbox{for } c=\xi_1+\xi_2\\[2mm]
\ds \frac{\tilde{r} }{m+\tilde{r}/r_j}+O\Bigr(\frac{\ep}{m+1}\Bigr)\quad &\mbox{for } c=\xi_j \ ( j=1,2).
\end{cases}
\eeq
Applying the mean-value property, we have
\beq\notag
\quad H\Bigr(0,0,\frac{a_\ep}{\xi_m^c}\Bigr)=
\begin{cases}
\ds H\Bigr(0,0,\frac{\tilde{r}}{m+1}\Bigr)+\frac{1}{m+1}O(\ep)\quad &\mbox{for } c=\xi_1+\xi_2\\[2mm]
\ds  H\Bigr(0,0,\frac{\tilde{r} }{m+\tilde{r}/r_j}\Bigr)+\frac{1}{m+1}O(\ep)\quad &\mbox{for }  c=\xi_j \ ( j=1,2),
\end{cases}
\eeq
Using Lemma \ref{CHC}, we have
\begin{align*}
&\ds\sum_{m=0}^\infty q_m^{\xi_1+\xi_2} H\left(\pm\mathbf{p}_m^{\xi_1+\xi_2}\right)
=\sum_{m=0}^\infty
\frac{4\pi\tilde{r}}{m+1}H\Bigr(0,0,\pm\frac{\tilde{r} }{m+1}\Bigr)+O({\ep} |\ln {\ep}|),\\
&\ds  \sum_{m=0}^\infty q_m^{\xi_j} H\left(\pm\mathbf{p}_m^{\xi_j}\right)=\sum_{m=0}^\infty
\frac{ 4\pi\tilde{r}}{m+\tilde{r}/r_j} H\Bigr(0,0,\pm\frac{\tilde{r} }{m+\tilde{r}/r_j} \Bigr)+O({\ep} |\ln {\ep}|), \quad j=1,2.
\end{align*}
Thanks to \eqnref{CHexp}, \eqnref{h_const} and the fact $\mu_1+\mu_2=1$, we prove \eqnref{eqn:CH}.

We can assume $H(\mathbf{0})=0$ as in the proof of Lemma \ref{CHC}, i.e.,
$$H(0,0,x_3)=\sum_{k=1}^\infty b_{H,k}x_3^k.$$
To make the notation simple, let us denote
$$z_m =  \frac{\tilde{r}}{m+1}\quad\mbox{and} \quad z_{m,j}=\frac{\tilde{r} }{m+\tilde{r}/r_j},$$
then \eqnref{CHdef} becomes
\begin{align}
\ds {C}_H &=
{4\pi\mu_1}\sum_{m=0}^\infty \sum_{k=1}^\infty b_{H,k}\left(z_m^{k+1}+(-z_{m,1})^{k+1}\right)
+{4\pi\mu_2}\sum_{m=0}^\infty \sum_{k=1}^\infty b_{H,k}\left(z_{m,2}^{k+1}+(-z_{m})^{k+1}\right)\notag\\
&={4\pi}\sum_{k=1}^\infty b_{H,k}\sum_{m=0}^\infty \Bigr[\mu_1z_m^{k+1}+\mu_2(-z_m)^{k+1}+
\mu_1(-z_{m,1})^{k+1}+\mu_2 z_{m,2}^{k+1} \Bigr].\label{CHdef2}
\end{align}
{We can exchange the order of summation in the above equation because of the absolute convergence of the series.}
From \eqnref{polygamma}, we immediately see that
\begin{align}
\sum_{m=0}^\infty z_m^{k+1} = \tilde{r}^{k+1} \zeta(k+1),
\qquad
\sum_{m=0}^\infty z_{m,j}^{k+1}= \tilde{r}^{k+1} \frac{(-1)^{k+1}}{k!}\psi_k(\tilde{r}/r_j).
\label{z_m}
\end{align}
Thus, from \eqnref{CHdef2} and \eqnref{z_m}, we prove (b).
\qed

\section{Asymptotics of $u$ and $\nabla u$}

In this section, we first approximate the series $S$ in \eqnref{def_S} by an integral defined in terms of bispherical coordinates, and then derive asymptotics for $h$, $u$, and $\nabla u$ in the Cartesian coordinates.

\subsection{Approximation of $h$ and $\nabla h$ by integrals in bispherical coordinates}\label{section:approximation}

Let us denote with $w_\theta$ in \eqnref{def:w_theta} that
\beq\label{def:g12}
\ds
\begin{cases}
\ds g(\xi,\theta;s)=\frac{w_\theta(\xi)}{w_\theta(\xi-s)},\\[3mm]
\ds g_1(\xi,\theta;s)=- \frac{\sinh(\xi-s)}{2a_\ep} \left(g(\xi,\theta;s)\right)^3,\\[3mm]
\ds  g_2(\xi,\theta;s)=\frac{\sinh{\xi}}{2a_\ep}g(\xi,\theta;s).
 \end{cases}
\eeq
Then we define
\begin{align}
\label{def_I}
\ds\mathcal{I}(\xi,\theta;c):&=\frac{1}{2(\xi_1+\xi_2)}\int_0^\infty g(\xi,\theta;t+c)\;dt, \\[2mm]
\ds\mathcal{I}_j(\xi,\theta;c):&=\frac{1}{2(\xi_1+\xi_2)}\int_0^\infty g_j(\xi,\theta;t+c)\;dt, \ j=1,2.\notag
\end{align}
It can be easily shown that
\begin{align}
\label{gjIj}\ds\Ical_j(\xi,\theta;c)&=
\begin{cases}\ds \frac{1}{2a_\ep(\xi_1+\xi_2)}\frac{w_\theta^3(\xi)}{w_\theta(\xi-c)}&\mbox{for }j=1\\[4mm]
\ds \frac{\sinh\xi}{2a_\ep}\Ical(\xi,\theta;c)&\mbox{for }j=2.
\end{cases}
\end{align}

Straightforward computations give
\begin{align}
\ds\notag\hat{\mathbf{e}}_\theta \cdot \nabla \big(g(\xi,\theta;s)\big)&=\frac{1}{\sigma_\theta(\xi,\theta)}\frac{\p}{\p\theta}\big(g(\xi,\theta;s)\big)\notag
=\frac{\sin{\theta}}{2a_\ep}\frac{w_\theta(\xi)\left(w^2_\theta(\xi-s)-w^2_\theta(\xi)\right)}{w^3_\theta(\xi-s)},\\
\label{eqn1015}\hat{\mathbf{e}}_\xi \cdot \nabla\big(g(\xi,\theta;s)\big)&=\frac{1}{\sigma_\xi(\xi,\theta)}\frac{\p}{\p\xi}\big(g(\xi,\theta;s)\big)=g_1(\xi,\theta;s)+g_2(\xi,\theta;s).\end{align}
Hence the directional derivatives of $\Ical$ at $(\xi,\theta,\varphi)$ becomes
\begin{align}
\ds\bigr(\hat{\mathbf{e}}_\theta \cdot\nabla \Ical\bigr)(\xi,\theta;c)
&=\frac{\sin\theta w_\theta(\xi)}{4a_\ep(\xi_1+\xi_2)}\int_0^\infty \frac{w^2_\theta(\xi-t-c)-w^2_\theta(\xi)}{w^3_\theta(\xi-t-c)}~ dt,\label{I_theta0}\\[1mm]
\label{eqn:exiIcal}
\bigr(\hat{\mathbf{e}}_\xi \cdot\nabla \Ical\bigr)(\xi,\theta;c)&=\Ical_1(\xi,\theta;s)+\Ical_2(\xi,\theta;s).
\end{align}
We also have
\begin{align}
\hat{\mathbf{e}}_\xi \cdot\nabla\bigr( \Ical(-\xi,\theta;s)\bigr)=\frac{1}{\sigma_\xi}\pd{}{\xi}\bigr( \Ical(-\xi,\theta;s)\bigr)
=-\Ical_1(-\xi,\theta;s)-\Ical_2(-\xi,\theta;s).\label{I_theta_m}
\end{align}

\begin{lemma}\label{nor_bdd}
There is a constant $C$ independent of $\ep$ such that
\begin{align}
\ds   \Bigr|\left(S-\Ical\right)\Bigr|(\xi,\theta;c)&\leq C\label{eqn:SIcal},\\
\ds \Bigr|\hat{\mathbf{e}}_\xi \cdot\nabla\left(S-\Ical\right)\Bigr|(\xi,\theta;c)&\leq C
\end{align}
for all $\theta\in[0,\pi]$ and $(\xi,c)\in[-\xi_1,\xi_2]\times\RR$ satisfying $|\xi|\leq(c-\xi)\leq 3(\xi_1+\xi_2)$.
\end{lemma}
\pf
Note that
$$
S(\xi,\theta;c)=\sum_{m=0}^\infty g(\xi,\theta;2 s_0 m +c), \quad s_0=\xi_1+\xi_2.
$$
By taking the directional derivative to the above and using \eqnref{eqn1015}, we obtain
\beq\label{S_deri_xi}
\bigr(\hat{\mathbf{e}}_\xi \cdot \nabla
S\bigr)(\xi,\theta;c) =
\sum_{m=0}^\infty g_1\big(\xi,\theta;2s_0m+c\big)+
\sum_{m=0}^\infty g_2\big(\xi,\theta;2s_0m+c\big).
\eeq
Since $ g_2(\xi,\theta;t+c) $ is positive {and} decreasing in $t\geq0$ owing to $(c-\xi)\geq |\xi|$, we derive
\beq\label{sumS1}
\bigg| \sum_{m=0}^\infty g_2\bigr(\xi,\theta;2s_0m+c\bigr)-\Ical_2(\xi,\theta,c)\bigg|\leq g_2(\xi,\theta;c) \leq 1.
\eeq
One can derive \eqnref{eqn:SIcal} by the same way.

To deal with the summation of $g_1(\xi,\theta;t+c)$ values in \eqnref{S_deri_xi}, which is not monotone in $t$, we now apply the Euler-{Mac}laurin summation formula, see section \ref{section:euler}. From Lemma \ref{app_dg2}, we estimate
\begin{align}
&\bigg| \sum_{m=0}^\infty g_1\bigr(\xi,\theta;2s_0m+c\bigr)-\Ical_1(\xi,\theta,c)\bigg|\nonumber\\
&\leq C \bigg(\big| g_1(\xi,\theta;c)\big|+{s_0}\bigg|\frac{\p  }{\p t} \bigr[g_1(\xi,\theta;t+c)\bigr]\Big|_{t=0}\bigg|+{s_0}\int_0^\infty \bigg|\frac{\p^2 }{\p t^2}\bigr[g_1(\xi,\theta;t+c)\bigr]\bigg|dt\bigg)\leq C\label{sumS2}
\end{align}
for some $C$ independent of $\ep$ and $(\xi,\theta)$.
{From \eqnref{eqn:exiIcal},\eqnref{S_deri_xi},\eqnref{sumS1}, and \eqnref{sumS2},
we prove the lemma.}\qed

\begin{lemma}\label{tan_KK_bdd}
There is a constant $C$ independent of $\epsilon$ such that
\begin{align*}
\ds\left|\hat{\mathbf{e}}_\theta \cdot\Bigr(\nabla \Ical\big(\xi,\theta;2(\xi_1+\xi_2)\big)-\nabla \Ical\big(-\xi,\theta;2\xi_1\big)\Big)\right|&\leq C,
\\[2mm]
\ds\left|\hat{\mathbf{e}}_\theta \cdot\Bigr(\nabla \Ical\big(-\xi,\theta;2(\xi_1+\xi_2)\big)-\nabla \Ical\big(\xi,\theta;2\xi_2\big)\Bigr)\right|&\leq C\quad\mbox{in }\RR^3\setminus\overline{(B_1\cup B_2)}.
\end{align*}
\end{lemma}
\pf
Applying \eqnref{I_theta0} and the mean value property, we have
\begin{align*}
\ds J&:=\left|\hat{\mathbf{e}}_\theta \cdot\Bigr(\nabla \Ical\big(\xi,\theta;2(\xi_1+\xi_2)\big)-\nabla \Ical\big(-\xi,\theta;2\xi_1\big)\Big)\right|\\
 \ds&=\bigg|\frac{ \sin\theta w_\theta(\xi) }{4a_\ep(\xi_1+\xi_2)}\int_{2\xi_1+\xi}^{2(\xi_1+\xi_2)-\xi} \frac{w^2_\theta(t)-w^2_\theta(\xi)}{ w^3_\theta(t)}dt\bigg|\\\ds& \leq C \bigg|  \frac{ \sin\theta w_\theta(\xi)}{a_\ep}\frac{w^2_\theta(\xi_0)-w^2_\theta(\xi)}{w^3_\theta(\xi_0)}\bigg|
\end{align*}
for some $\xi_0 \in \bigr(2\xi_1+\xi,\;2(\xi_1+\xi_2)-\xi\bigr)$ and a constant $C$ independent of $\ep$ and $(\xi,\theta)$.
By applying the mean value property again, we have
$$
w_\theta(\xi_0)-w_\theta(\xi)=(\xi_0-\xi) \frac{\sinh\xi_*}{2 w_\theta(\xi_*)}\quad\mbox{for some }\xi_*\in(\xi,\xi_0).
$$
Note that $|\xi|\leq|\xi_0|$ and, hence,
$$
{w_\theta(\xi)}\leq{w_\theta(\xi_0)}. $$
Therefore, we conclude
\begin{align*}
\ds J &\leq C\frac{ |\sin\theta||\xi_0-\xi|\sinh\xi_*}{a_\ep w_\theta(\xi_*)}
\left( \frac{1}{w_\theta(\xi_0)}+\frac{w_\theta(\xi)}{w_\theta(\xi_0)^2}\right)\leq C \frac{|\sin\theta|\sinh \xi_*}{w_\theta(\xi_0) w_\theta(\xi_*)}\leq C.\end{align*}

 Similarly, we can prove the second uniform boundedness.
\qed

\subsection{Asymptotics of $h$ and $u$ in the bispherical coordinates}

Since $h$ is a linear combination of $S$, see Lemma \ref{h_L}, a direct consequence of the previous lemmas is {the asymptotics} of $h$ and $\nabla h$ in terms of integrals. We fix some notations for the sake of notational simplicity before {deriving the asymptotics}: let us denote
$$\tilde{\mu}_j=\mu_\ep\hskip .3mm\mu_j,\ j=1,2,$$
and \begin{align}
\ds h_s(\Bx) &= -\tilde\mu_1\Bigr(\mathcal{I}(-\xi,\theta;2\xi_1)-\mathcal{I}(\xi,\theta;2\xi_1+2\xi_2)\Bigr)+\tilde\mu_2\Bigr( \mathcal{I}(\xi,\theta;2\xi_2)- \mathcal{I}(-\xi,\theta;2\xi_1+2\xi_2)\Bigr),\notag\\
q_h(\Bx)&=\tilde\mu_1\Bigr( \mathcal{I}_1(-\xi,\theta;2\xi_1)+\mathcal{I}_1(\xi,\theta;2\xi_1+2\xi_2)\Bigr)+\tilde\mu_2 \Big(\mathcal{I}_1(\xi,\theta;2\xi_2)+ \mathcal{I}_1(-\xi,\theta;2\xi_1+2\xi_2)\Bigr).\label{qhint}
\end{align}

\begin{prop}\label{hh0_09081}
The solution $h$ to {\rm\eqnref{h:eqn} }satisfies
\begin{align}
\ds h(\Bx)&=h_s(\Bx)+b(\Bx),\label{eqn:1013_1}\\[1mm]
\ds\nabla h(\Bx) &= q_h(\Bx)\hat{\mathbf{e}}_\xi(\Bx)+r(\Bx),\label{eqn:1013_2}\end{align}
where $\|b\|_\infty$, $\|\nabla b\|_\infty$ and $\|r\|_\infty$ are bounded independently of $\ep$.
Moreover, we have
\beq\label{eqn:Qbound}
C_1|\ep\ln\ep|^{-1}\leq \|q_h\|_\infty\leq C_2|\ep\ln\ep|^{-1},
\eeq
for some positive constants $C_1$ and $C_2$ independent of $\ep$.
\end{prop}
\pf
Firstly, we prove that $\|\nabla b \|_\infty=\|\nabla (h - {h}_s)\|_\infty $ is uniformly bounded regardless of $\ep>0$.
From Lemma \ref{def_g2}, $ h_s$ is harmonic and has the decay property at infinity. Hence it is enough to derive the uniform boundedness of $|\nabla(h- h_s)|$ in $\Bx\in\p B_1 \cup \p B_2$ and $\ep>0$. For simplicity, we consider only for $\p B_1$. Since $h$ is constant on $\p B_1$, we have
\begin{align}\label{eqn:hminushs}
\big|\nabla(h- h_s)\big| \leq  \big|\hat{\mathbf{e}}_\xi \cdot \nabla(h-  h_s)\big|+\big|\hat{\mathbf{e}}_\theta \cdot \nabla  h_s\big| \quad\mbox{on }\p B_1.
\end{align}
In the following we show that $\|\hat{\mathbf{e}}_\xi \cdot \nabla(h- h_s)\|_\infty$ and $\|\hat{\mathbf{e}}_\theta \cdot \nabla  h_s \|_\infty$ are uniformly bounded in $\ep$.
Note that the directional derivatives of $h$ and $h_s$ are combinations of those of $\Ical$. More precisely speaking, because of $\hat{\mathbf{e}}_\xi \cdot\nabla \bigr(S(-\xi,\theta;s)\bigr)=-\hat{\mathbf{e}}_\xi \cdot\bigr(\nabla S(-\xi,\theta;s)\bigr)$ and $\hat{\mathbf{e}}_\xi \cdot\nabla \bigr(\Ical(-\xi,\theta;s)\bigr)=-\hat{\mathbf{e}}_\xi \cdot\bigr(\nabla \Ical(-\xi,\theta;s)\bigr)$, one can rewrite $\hat{\mathbf{e}}_\xi \cdot\nabla h$ and $\hat{\mathbf{e}}_\xi \cdot\nabla h_s$ as
\begin{align}
\ds\hat{\mathbf{e}}_\xi\cdot\nabla  h =\;&-C_1 \hat{\mathbf{e}}_\xi \cdot\Bigr(\nabla S\big(-\xi,\theta;2\xi_1\big)+\nabla S\big(\xi,\theta;2\xi_1+2\xi_2\big)\Bigr)\notag \\\ds&\ds +C_2 \hat{\mathbf{e}}_\xi \cdot\Bigr(\nabla S\big(\xi,\theta;2\xi_2\big)+\nabla S\big(-\xi,\theta;2\xi_1+2\xi_2\big)\Big), \label{xi_h}\\
\ds\hat{\mathbf{e}}_\xi\cdot\nabla  h_s=\;&\ds\tilde\mu_1 \hat{\mathbf{e}}_\xi \cdot\Bigr(\nabla \Ical\big(-\xi,\theta;2\xi_1\big)+\nabla \Ical\big(\xi,\theta;2\xi_1+2\xi_2\big)\Bigr)\notag \\\ds+&\ds \tilde\mu_2\hat{\mathbf{e}}_\xi \cdot\Bigr(\nabla \Ical\big(\xi,\theta;2\xi_2\big)+\nabla \Ical\big(-\xi,\theta;2\xi_1+2\xi_2\big)\Big). \label{xi_hst}
\end{align}
Similarly to \eqnref{xi_hst}, we use the fact
$\hat{\mathbf{e}}_\theta \cdot\nabla \bigr(\Ical(-\xi,\theta;s)\bigr)=\hat{\mathbf{e}}_\theta \cdot\bigr(\nabla \Ical(-\xi,\theta;s)\bigr)$ to rewrite $\hat{\mathbf{e}}_\theta\cdot\nabla h_s$ as
\begin{align}
\ds\hat{\mathbf{e}}_\theta\cdot\nabla h_s=&-\tilde{\mu}_1\hat{\mathbf{e}}_\theta \cdot\Bigr(\nabla \Ical\big(-\xi,\theta;2\xi_1\big)-\nabla \Ical\big(\xi,\theta;2\xi_1+2\xi_2\big)\Bigr)\notag\\\ds+&\tilde{\mu}_2\hat{\mathbf{e}}_\theta \cdot\Bigr(\nabla \Ical\big(\xi,\theta;2\xi_2\big)-\nabla \Ical\big(-\xi,\theta;2\xi_1+2\xi_2\big)\Big).\label{theta_hst}
\end{align}

Suppose that $(\tilde\xi,c)$ is one of $(-\xi,2\xi_1),(\xi,2\xi_2),(\pm\xi,2\xi_1+2\xi_2)$ with $-\xi_1\leq\xi\leq\xi_2$. Then we have \beq\notag|\tilde\xi|\leq (c-\tilde\xi)\leq3(\xi_1+\xi_2).\eeq
Since $|\tilde\xi-c|\geq|\tilde\xi|$, we have $0<g(\tilde\xi,\theta,c)\leq 1$. Using this and the definition of $\Ical$ and $\Ical_j$'s, we can easily show
\beq \label{I1I2bound}0< \Ical_1(\tilde\xi,\theta;c)\leq \frac{C}{\ep},
\quad 0<  \Ical(\tilde\xi,\theta;c), \Ical_2(\tilde\xi,\theta;c) \leq \frac{C}{\sqrt\ep}\quad\mbox{in }\RR^3\setminus(B_1\cup B_2).\eeq
Here and in the remaining of the proof, $C$ indicates a positive constant independent of $\ep$ and $(\xi,\theta)$.
Thanks to \eqnref{eqn:exiIcal}, one obtains
\beq\label{ineq0908}\hat{\mathbf{e}}_\xi \cdot\nabla \Ical(\tilde\xi,\theta;c)=O(\ep^{-1}).\eeq
We also have from Lemma \ref{nor_bdd} that
\beq \label{nor_bdd_2}
\Big| \hat{\mathbf{e}}_\xi \cdot\nabla (S-\Ical) \Big|(\tilde\xi,\theta;c) \leq C.\eeq
Note that
\beq\notag
\tilde\mu_1,\tilde\mu_2=O(|\ln\ep|^{-1})\quad\mbox{and}\quad\tilde\mu_j=(-1)^jC_j+O(\ep).\eeq
Using these facts, \eqnref{ineq0908} and \eqnref{nor_bdd_2}, we get
\begin{align}
&\Big| \hat{\mathbf{e}}_\xi \cdot \left(C_j \nabla S
- (-1)^j\mu_j \nabla \Ical\right) \Big|(\tilde\xi,\theta;c) \notag\\&\leq
\Big| C_j \hat{\mathbf{e}}_\xi \cdot  \nabla \left(S
-   \Ical\right)\Big|(\tilde\xi,\theta;c) +\Big|\left(C_j-(-1)^j\tilde\mu_j\right) \hat{\mathbf{e}}_\xi \cdot \nabla \Ical\Big|(\tilde\xi,\theta;c) \notag\\&\leq C.\notag
\end{align}
Hence we obtain from \eqnref{xi_h} and \eqnref{xi_hst} that
\beq\label{eqn:nablahhs}\big\|\hat{\mathbf{e}}_\xi \cdot \nabla(h- h_s)\big\|_\infty \leq C.\eeq
The $\theta$-directional derivative of $h_s$ satisfies
\beq\label{hs_tan_bdd}
\big\|\hat{\mathbf{e}}_\theta \cdot \nabla  h_s \big\|_\infty \leq C
\eeq
due to \eqnref{theta_hst} and Lemma \ref{tan_KK_bdd}.
Thanks to \eqnref{eqn:hminushs}, \eqnref{eqn:nablahhs} and \eqnref{hs_tan_bdd}, we derive that that $\|\nabla(h- h_s)\|_\infty$ is uniformly bounded independently of $\ep$. This shows that $\|\nabla b\|_\infty\leq C$ by the discussion at the beginning of the proof.
In fact, due to \eqnref{hs_tan_bdd}, we have shown a slightly stronger result as follows:
\beq\label{h_hsxi}
\nabla h(\Bx) = (\hat{\mathbf{e}}_\xi \cdot \nabla h_s)(\Bx) \hat{\mathbf{e}}_\xi(\Bx) +
\tilde{r}(\Bx),\quad \| \tilde{r}\|_\infty\leq C.
\eeq

Now we prove \eqnref{eqn:1013_2}. From \eqnref{eqn:exiIcal} and the definition of $\Ical_j$'s, the $\xi$-directional derivative of $h_s$ satisfies
\beq\label{hsxi_q_v}
\hat{\mathbf{e}}_\xi\cdot\nabla h_s(\Bx)=q_h(\Bx)+v(\Bx),
\eeq
where $q_h$ is defined as in \eqnref{qhint} and
$$v(\Bx)=\frac{\sinh \xi}{2a_\ep}{h}_s(\Bx).$$
We need to show that $\| v\|_\infty$ is bounded regardless of $\ep>0$. To do that let us consider the remainder term $b$ in \eqnref{eqn:1013_1}.
From \eqnref{eqn:SIcal} and Lemma \ref{h_L}, one can easily prove
\beq\label{h_hs_bdd}
\|b\|_\infty=\|h-h_s\|_\infty\leq C
\eeq
similarly to the proof of \eqnref{eqn:nablahhs}.
Remind that $h$ has the decaying property and $h\bigr|_{\p B_1\cup \p B_2}=O(|\ln\ep|^{-1})$. Hence, $\|h\|_\infty$ is bounded independently of $\ep$ and so does for $\|h_s\|_\infty$ thanks to \eqnref{h_hs_bdd}. So we have
\beq\notag
\|v\|_\infty=\Big\| \frac{\sinh \xi}{2a_\ep}{h}_s\Big\|_\infty\leq C.
\eeq
Hence, we obtain \eqnref{eqn:1013_2} using \eqnref{h_hsxi} and \eqnref{hsxi_q_v}.

We note from the definition of $q_h$ that $q_h(\Bx_0)$ for $\Bx_0\in\p B_1$ of which bispherical coordinates are $(-\xi_1,\pi,0)$ satisfies
\begin{align*}
q_h(\Bx_0)\geq \tilde{\mu}_1 \Ical_1(\xi_1,\pi;2\xi_1)=\frac{\tilde{\mu}_1}{2a_\ep(\xi_1+\xi_2)}\frac{(\cosh\xi_1-\cos\pi)^{\frac{3}{2}}}{(\cosh\xi_1-\cos\pi)^{\frac{1}{2}}}\geq \frac{C}{\ep|\ln\ep|}.
\end{align*}
This proves the lower bound in \eqnref{eqn:Qbound}, and the upper bound follows from \eqnref{I1I2bound}. Hence we finish the proof.
\qed

We now have the asymptotics of $u$ and $\nabla u$ thanks to
\eqnref{bdd_g_def}, \eqnref{eqn:CH} and Proposition \ref{hh0_09081} as follows.
\begin{prop}\label{uasymp:bispherical}
The solution $u$ to {\rm\eqnref{u:eqn}} satisfies
\begin{align*}
\ds u(\Bx)&=  C_H h_s(\Bx) + H(\Bx)+b(\Bx),\\[1mm]
\ds\nabla u(\Bx)&=C_H q_h(\Bx)\hat{\mathbf{e}}_\xi(\Bx)+\nabla H(\Bx)+r(\Bx),
\end{align*}
where $\|b\|_\infty$, $\|\nabla b\|_\infty$ and $\|r\|_\infty$ are bounded independently of $\ep$.
\end{prop}

\subsection{Asymptotics of $h$ and $u$ in the Cartesian coordinates}

With the notations defined in \eqnref{def:tilder} and \eqnref{mu1002}, we define two density functions $\rho_j$, $j=1,2$, as
\begin{align}
\ds \rho_1 (0,0,c)&=\frac{\tilde r {\mu_\ep} }{\sqrt{c^2-a_\ep^2}}\Bigr(\mu_1\mathbbm{1}_{[\Bc_1,\Bp_1]}+\mu_2\mathbbm{1}_{[R_1(\Bc_2),\Bp_1]}\Bigr)(0,0,c),\\
\ds \rho_2 (0,0,c)&=\frac{\tilde r {\mu_\ep} }{\sqrt{c^2-a_\ep^2}}\Bigr(\mu_2\mathbbm{1}_{[\Bp_2,\Bc_2]}+\mu_1\mathbbm{1}_{[\Bp_2,R_2(\Bc_1)]}\Bigr)(0,0,c),
\end{align}
where the symbol $[\Bx_1,\Bx_2]$ means the line segment connecting two points $\Bx_1$ and $\Bx_2$, and $\mathbbm{1}_{[\Bx_1,\Bx_2]}$ is the indicator function of $[\Bx_1,\Bx_2]$.  From Proposition \ref{uasymp:bispherical}, we derive the following corollary which tells that the solution $h$ to \eqnref{h:eqn} can be expressed as the integral with the integrand $\rho_1$ and $\rho_2$.

\begin{cor}\label{lem:hhs}
We have
$$ h(\Bx)= -\int_{[\Bc_1,\Bp_1]}\frac{\rho_1(\Bc)}{|\Bx-\Bc|}d\Bc
+\int_{[\Bp_2,\Bc_2]}\frac{\rho_2(\Bc)}{|\Bx-\Bc|}d\Bc+r(\Bx)\quad\mbox{in }\RR^3\setminus\overline{B_1\cup B_2},$$
where $\|\nabla r\|_\infty$ is bounded regardless of $\ep$.
\end{cor}

\pf
Let us express the function $h_s$ {in} the Cartesian coordinates.
Applying Lemma \ref{def_g2} and letting $c=a_\ep\coth(t+\xi_1)$, one computes
\begin{align}
\Ical(-\xi,\theta;2\xi_1)\notag
&=\frac{1}{\xi_1+\xi_2}\int_0^\infty\frac{w_\theta(\xi)}{w_\theta(\xi+2t+2\xi_1)}~dt\\\notag
&=\frac{a_\ep}{\xi_1+\xi_2}\int_0^\infty\frac{1}{\sinh(t+\xi_1)}\frac{1}{|\Bx+a_\ep\coth(t+\xi_1)\mathbf{e}_3|}~dt\\\notag
&=\frac{a_\ep}{\xi_1+\xi_2}\int_{a_\ep\coth \xi_1}^{a_\ep}  \frac{-1}{\sqrt{c^2-a_\ep^2}}\frac{1}{|\Bx+c\mathbf{e}_3|}~dc\\
&=\int_{[\Bc_1,\Bp_1]}f(\Br)d\Br,\notag
\end{align}
where $$f(\Br) = \frac{a_\ep}{\xi_1+\xi_2}\frac{1}{\sqrt{|\Br|^2-a_\ep^2}}\frac{1}{|\Bx-\Br|}.$$
Similarly, one can easily obtain
\begin{align}
\ds\Ical(\xi,\theta;2\xi_2) &= \int_{[\Bp_2,\Bc_2]}f(\Br)d\Br,\notag\\
\ds\Ical(-\xi,\theta;2(\xi_1+\xi_2))& = \int_{[R_1(\Bc_2),\Bp_1]}f(\Br)d\Br,\notag\\
\ds\Ical(\xi,\theta;2(\xi_1+\xi_2)) &= \int_{[\Bp_2,R_2(\Bc_1)]}f(\Br)d\Br.\notag
\end{align}
Hence we have
\beq\label{hs09091}h_s(\Bx)=\frac{a_\ep}{\tilde r(\xi_1+\xi_2)}\tilde{h}_s(\Bx),\eeq
where $$\tilde{h}_s(\Bx)= -\int_{[\Bc_1,\Bp_1]}\frac{\rho_1(\Bc)}{|\Bx-\Bc|}d\Bc
+\int_{[\Bp_2,\Bc_2]}\frac{\rho_2(\Bc)}{|\Bx-\Bc|}d\Bc.
$$
Thanks to Proposition \ref{hh0_09081} and the fact $\frac{a_\ep}{\tilde r(\xi_1+\xi_2)}=1+O(\ep)$, we prove the corollary.
 \qed

The following corollary is the direct consequence of \eqnref{bdd_g_def}, Theorem \ref{lem:CH} and Corollary \ref{lem:hhs}.
\begin{cor}
The solution $u$ to {\rm\eqnref{u:eqn}} satisfies
$$u(\Bx)=  C_H \left( -\int_{[\Bc_1,\Bp_1]}\frac{\rho_1(\Bc)}{|\Bx-\Bc|}d\Bc
+\int_{[\Bp_2,\Bc_2]}\frac{\rho_2(\Bc)}{|\Bx-\Bc|}d\Bc\right) + H(\Bx)+b(\Bx),$$
where $\|\nabla b\|_\infty$ is bounded independently of $\ep$.
\end{cor}

 Near the fixed points $\Bp_1$ and $\Bp_2$, {the density functions $\rho_1$ and $\rho_2$ are of similar form} to that obtained in \cite{Pol} as mentioned in the introduction.
It is worth to emphasize that, in this paper, we derived the continuous image charge distribution by rigorous asymptotic analysis without any physical assumptions.
Moreover, it turns out that each of density functions $\rho_1$ and $\rho_2$ has an discontinuity, and the coefficients in the density functions are explicitly calculated.

\subsection{Proof of Theorem \ref{main_thm2}}\label{subset:proof theorem}

\noindent{{\textbf{Proof of Theorem \ref{main_thm2}}}}
From \eqnref{partialBj}, \eqnref{xi0909}, \eqnref{exi:x} and Lemma \ref{def_g2}, we have
\begin{align*}\Ical_1(-\xi,\theta;2\xi_1)\hat{\mathbf{e}}_\xi
&=\frac{1}{2a_\ep(\xi_1+\xi_2)}\frac{w_\theta^3(\xi)}{w_\theta(\xi+2\xi_1)}\hat{\mathbf{e}}_\xi
=\frac{\cosh \xi-\cos\theta}{2a_\ep(\xi_1+\xi_2)}\frac{w_\theta(\xi)}{w_\theta(\xi+2\xi_1)}\frac{a_\ep}{\cosh\xi-\cos\theta}\mathbf{N}(\Bx)\\
&=\frac{1}{2(\xi_1+\xi_2)}\frac{a_\ep}{\sinh\xi_1}\frac{\mathbf{N}(\Bx)}{|\Bx-\Bc_1|}\\
&=\frac{r_1}{2(\xi_1+\xi_2)}\frac{\mathbf{N}(\Bx)}{|\Bx-\Bc_1|}
\end{align*}
and, by the same way,
$$\Ical_1(\xi,\theta;2\xi_2)\hat{\mathbf{e}}_\xi=\frac{r_2}{2(\xi_1+\xi_2)}\frac{\mathbf{N}(\Bx)}{|\Bx-\Bc_2|}.$$
Similarly, we compute
\begin{align*}
\Ical_1(-\xi,\theta;2\xi_1+2\xi_2)\hat{\mathbf{e}}_\xi&=\frac{1}{2(\xi_1+\xi_2)}\frac{a_\ep}{\sinh(\xi_1+\xi_2)}\frac{\mathbf{N}(\Bx)}{|\Bx-R_1(\Bc_2)|}\\
&=\frac{\tilde r}{2(\xi_1+\xi_2)}\frac{\mathbf{N}(\Bx)}{|\Bx-R_1(\Bc_2)|}\Bigr(1+O(\ep)\Bigr),\\
\Ical_1(+\xi,\theta;2\xi_1+2\xi_2)\hat{\mathbf{e}}_\xi&=\frac{\tilde r}{2(\xi_1+\xi_2)}\frac{\mathbf{N}(\Bx)}{|\Bx-R_2(\Bc_1)|} \Bigr(1+O(\ep)\Bigr).
\end{align*}
Therefore we have
\beq\notag
q_h(\Bx)\hat{\mathbf{e}}_\xi(\Bx)= \frac{a_\ep}{\tilde r(\xi_1+\xi_2)}\psi(\Bx)\mathbf{N}(\Bx)+O(1).
 \eeq
Note that $\frac{\tilde r(\xi_1+\xi_2)}{a_\ep}=1+O(\ep).$
From Proposition \ref{uasymp:bispherical} and \eqnref{qhint}, we prove Theorem \ref{main_thm2}.
\qed

\smallskip

\noindent{{\textbf{Proof of Theorem \ref{prop:superfocusing}}}}
To prove (a), it is enough to show, in view of Proposition \ref{hh0_09081}, that there is a constant $C$ independent of $\ep$ such that
$$
| q_h(\xi,\theta)| \leq C\quad \mbox{for } |\theta|\leq\sqrt{\ep|\ln\ep|},\;-\xi_1\leq\xi\leq\xi_2.
$$
 We estimate only one term $\mu_1 \mathcal{I}_1(-\xi,\theta;2\xi_1)$ in $h_1$; the other three terms can be estimated in a similar way.
The bispherical {coordinates} $(\xi,\theta,\varphi)$ of $\Bx\in \mathbb{R}^3\setminus\overline{(B_1\cup B_2\cup{\Omega_{\ep}^*})}$ satisfies $|\theta|\leq\sqrt{\ep|\ln\ep|}$ and $-\xi_1\leq\xi\leq\xi_2$, so that it follows
\begin{gather*}
 |\xi+2\xi_1|\geq|\xi|,\\
 w^2_\theta(\xi) =\cosh\xi-\cos\theta = 1+O(\ep)-\big(1+O(\ep|\ln\ep|)\big)=O(\ep|\ln\ep|).
\end{gather*}
We compute
$$
\big| \mu_\ep\mu_1 \mathcal{I}_1(-\xi,\theta;2\xi_1)\big|=\frac{\mu_\ep\mu_1}{2a_\ep(\xi_1+\xi_2)}\frac{w_\theta(\xi)}{w_\theta(\xi+2\xi_1)}w^2_\theta(\xi)\leq C
 \frac{\mu_\ep\mu_1}{a_\ep(\xi_1+\xi_2)} \ep|\ln\ep| \leq C,
$$
where $C$ is a constant independent of $\ep$. This proves (a).

\smallskip

We prove (b) by showing that $
|q_h(\xi,\tilde{\theta}_\ep)|\rightarrow \infty$ as $\ep$ tends to zero.
Again, we consider only $\mu_1\mathcal{I}_1(-\xi,\tilde{\theta}_\ep;2\xi_1)$.
Let us denote $w_\theta(\xi)=w(\xi,\theta)$ for notational sake. Because of $-\xi_1\leq \xi \leq \xi_2$,
we have $$
w^2(\xi,\tilde{\theta}_\ep)=\cosh\xi-\cos\tilde{\theta}_\ep=\frac{\tilde{\theta}_\ep^2}{2}+O\big(\tilde{\theta}_\ep^4\big)
$$
and
$$
\frac{w(\xi,\tilde{\theta}_\ep)}{w(\xi+2\xi_1,\tilde{\theta}_\ep)}
=\left(\frac{\cosh\xi-\cos\tilde{\theta}_\ep}{\cosh(\xi+2\xi_1)-\cos\tilde{\theta}_\ep}\right)^{\frac{1}{2}}
=\left( \frac{\tilde{\theta}_\ep^2+O(\tilde{\theta}_\ep^4)}{\tilde{\theta}_\ep^2+O(\tilde{\theta}_\ep^4)}\right)^{\frac{1}{2}}
\longrightarrow 1\quad \mbox{as }\ep\rightarrow 0.
$$
Hence, we have \begin{align*}
\big| \mu_\ep\mu_1 \mathcal{I}_1(-\xi,\tilde{\theta}_\ep;2\xi_1)\big|
&=
\left|\frac{\mu_\ep\mu_1}{2a_\ep (\xi_1+\xi_2)}\frac{w^3(\xi,\tilde{\theta}_\ep)}{w(\xi+2\xi_1,\tilde{\theta}_\ep)}\right|\geq C \frac{\tilde{\theta}_\ep^2}{\ep|\ln\ep|}\longrightarrow \infty\quad \mbox{as }\ep\rightarrow 0.
\end{align*}
Here, we can choose $C$ independent of $\ep$ and $\xi$ satisfying $-\xi_1\leq\xi\leq\xi_2$.
This proves (b).
\qed

\section{The Euler-Maclaurin formula and its two applications}\label{section:euler}
The following is a special case of the Euler-Maclaurin summation formula: for $f \in C^2[x_0,\infty)$ satisfying $f,f',f'' \in L^1(x_0,\infty)$, we have for any ${s}>0$ that
$$
\sum_{k=0}^{\infty} f(x_0+k \tilde{s}) \tilde{s} = \int_{x_0}^\infty f(x)dx +\frac{{s}}{2}f(x_0)-\frac{{s}^2}{12}f'(x_0) + {R},$$
where the remainder term satisfies
$$
|{R}| \leq \frac{{s}^2}{12} \int_{x_0}^\infty |f''(x)| dx.
$$

The followings are the applications of the Euler-Maclaurin formula, and they are essentially used to prove the main theorems in this paper.
\begin{lemma} \label{Sum_fk_lem}
Fix $s_0,c$ such that $0<s_0,c\leq (\xi_1+\xi_2)$ with $\xi_1,\xi_2$ given by {\rm\eqnref{xi0909}} and set
\beq\label{fk}f_k(x)=\frac{\coth^k x}{\sinh x}-\frac{1}{x^{k+1}} \quad \mbox{for }x>0,~k=1,2,\dots.\eeq
Then there is a constant $C$ independent of $\ep$ and $k$ such that
$$
\left|s_0\sum_{m=0}^\infty f_k(s_0 m+c) \right| \leq C 2^k \ep^{\frac{-k+2}{2}}|\ln \ep|
$$
for small enough $\ep>0$.
\end{lemma}
\pf
Let $k$ be a fixed positive integer. One can easily check
\begin{align*}
\ds
f_k'(x)&=\frac{k+1}{ x^{k+2}}-\frac{\coth^{k+1}x}{\sinh x}- \frac{k\coth^kx}{\sinh^2x \cosh x},\\
\ds f_k''(x)&=-\frac{(k+2)(k+1)}{x^{k+3}}+\frac{\coth^{k+2}x}{\sinh x}+\frac{(\coth^kx) \bigr[(4k+1)\cosh^2 x+k(k-1) \bigr]}{\sinh^3x\cosh^2 x}.
\end{align*}
Applying the Euler-Maclaurin summation formula, we have
\beq\notag\left|\sum_{m=0}^\infty f_k(c+s_0 m) s_0\right|\leq \left| \int_c^\infty f_k(x) dx \right| + \frac{s_0}{2}  |f_k(c)|+\frac{s_0^2}{12}|f_k'(c)|+\frac{s_0^2}{12}\int_c^\infty |f_k''(t)|dt.
\eeq
In the following, we estimate the four terms in the right-hand side in the equation above.

We first define a function $v_k$ as
$$
v_k(x)=\frac{1}{x^2}-\frac{\tanh^k x}{x^{k+2}} \quad \mbox{for }  x>0.
$$
Note that $v_k$ is positive and monotonically decreasing and
\beq\label{vkzero}\lim_{x\rightarrow 0^+} v_k(x)=\frac{k}{3}.\eeq
To estimate $|f_k(c)|$ and $\left| \int_c^\infty f_k(x) dx \right|$, we decompose $f_k$ for $k\in\NN$ as
\beq\notag
f_k(t)=\frac{1}{t^{k-1}}\frac{t^k}{\tanh^k t}\bigr(v_k(t) +p_0(t) \bigr)\quad\mbox{with }p_0(t)=\frac{1}{t \sinh t}-\frac{1}{t^2}.
\eeq
Note that $p_0(t)$ is bounded on $\{t>0\}$ and $({c^k}/{\tanh^k c})$ is bounded by a constant $C$. Here and in the remaining of the proof, $C$ indicates a constant independent of $\ep$ and $k$. Hence we have
\beq\label{fk00}
\ds|f_k(c)| \leq   \frac{C}{c^{k-1}}  \Bigr(|p_0(c)|+|v_k(c)|\Bigr)\leq C k \ep^{\frac{-k+1}{2}}.
\eeq

Let us now estimate $\left| \int_c^\infty f_k(x) dx \right|$. For $k=1$, we have \beq\notag
\left|\int_c^\infty f_1(x) dx\right| \leq \int_0^c \left|f_1(x)\right| dx \leq C \sqrt\ep
\eeq
thanks to the boundedness of $f_1(x)$ on $\{t>0\}$ and $\int_0^\infty f_1(x)dx =0$. For $k\geq2$, we have
\begin{align*}
\ds\left| \int_c^\infty f_k(x) dx \right|&\leq  \int_c^1 \left|f_k(x)\right|  dx + \int_1^\infty \left| f_k(x)\right| dx\notag \\
\ds&\leq \int_c^1 \frac{1}{x^{k-1}}\left(\frac{x}{\tanh x}\right)^k (|p_0(x)|+|v_k(x)|) dx+ C
\int_1^\infty \frac{1}{\sinh x} dx +\frac{1}{k}\nonumber \\
\ds& \leq C k \left(\frac{1}{\tanh 1}\right)^k \int_c^1 x^{-k+1} dx +C\leq C 2^k \ep^{\frac{-k+2}{2}}|\ln \ep|.
\end{align*}
Hence it follows
\begin{align}
\ds\left| \int_c^\infty f_k(x) dx \right|\leq C 2^k \ep^{\frac{-k+2}{2}}|\ln \ep|\quad\mbox{for }k\in\NN.\label{fk03}
\end{align}

Similar to the decomposition of $f_k$, we have such decompositions of $f_k'$ and $f_k''$:
\begin{align*}
\ds
f_k'(t)&=\frac{1}{t^{k}}\frac{t^k}{\tanh^k t}\Bigr(-(k+1)v_k(t)+p_1(t)+k p_2(t)\Bigr),\\
\ds f_k''(t)&=\frac{1}{t^{k+1}}\frac{t^k}{\tanh^k t}\Bigr( (k^2+3k+2)v_k(t)+p_3(t)+(4k+1)p_4(t)
+k(k-1)p_5(t)  \Bigr)
\end{align*}
with bounded functions
\begin{align*}
p_1(t)&= \frac{1}{t^2}-\frac{\coth t}{\sinh t}, \quad p_2(t)=\frac{1}{t^2}-\frac{1}{\sinh^2 t\cosh t},\\
p_3(t)&=\frac{t \coth^2 t}{\sinh t}-\frac{1}{t^2}, \quad p_4(t)=\frac{t}{\sinh^3 t}-\frac{1}{t^2}, \quad p_5(t)=\frac{t}{\cosh^2t \sinh^3 t}-\frac{1}{t^2}.
\end{align*}
Using \eqnref{vkzero}, we derive
\begin{align}\label{fk01}
\ds|f_k'(c)| &\leq C \frac{1}{c^{k}} \bigr((k+1)|v_k(c)|+|p_1(c)|+k|p_2(c)|\bigr) \leq C k^2 \epsilon^{-\frac{k}{2}}
\end{align}
and
\begin{align}
\ds \int_c^\infty |f_k''(t)| dt &= \int_c^1 |f_k''(t)|dt+\int^\infty_1 |f_k''(t)|dt\notag\\
&\leq C \int_c^1 \frac{1}{t^{k+1}} \frac{k^3}{\tanh^k 1}dt +C \int_1^\infty \left( \frac{k^2+3k+2}{t^{k+3}}+ \frac{1}{\sinh t}
+ \frac{k^2+3k+1}{\sinh^3 t}\right)dt\notag\\
&\leq  C \left(2^k \epsilon^{-\frac{k}{2}} +k^2\right)\leq  C 2^k \epsilon^{-\frac{k}{2}}\label{fk02}.
\end{align}

From \eqnref{fk00}, \eqnref{fk03}, \eqnref{fk01}, and \eqnref{fk02}, the conclusion follows.
\qed

\begin{lemma}\label{app_dg2}
Let $c$ be a constant satisfying $|\xi|\leq(c-\xi)\leq 3(\xi_1+\xi_2)$ for all $-\xi_1\leq \xi \leq \xi_2$ and
$g_1$ be the function given by {\rm\eqnref{def:g12}}, i.e., \beq \notag g_1(\xi,\theta;c)= -\frac{\sinh(\xi- c)}{2a_\ep}  \left(\frac{w_\theta(\xi)}{w_\theta(\xi-s)}\right)^3.\eeq
Then there exists a constant $C$ independent of $\xi$, $\theta$, $\ep$ and $c$ such that
\beq\notag
\bigg| \sum_{m=0}^\infty {g_1}\bigr(\xi,\theta;2s_0m+c\bigr)-\frac{1}{2s_0}\int_0^{\infty} {g_1}(\xi,\theta;t+c) dt\bigg| \leq C\quad\mbox{with }s_0=\xi_1+\xi_2,\eeq
for all $-\xi_1\leq \xi \leq \xi_2$ and $0\leq\theta\leq \pi$.
\end{lemma}

\pf
We use the Euler-{Mac}laurin summation formula to have \begin{align*}
&\left| \sum_{m=0}^\infty {g_1}\bigr(\xi,\theta;s_0m+c\bigr)-\frac{1}{s_0}\int_0^{\infty} {g_1}(\xi,\theta;t+c) dt\right|\nonumber\\
&\leq\big| {g_1}(\xi,\theta;c)\big|+s_0\bigg|\frac{\p  }{\p t} \bigr[{g_1}(\xi,\theta;t+c)\bigr]\Big|_{t=0}\bigg|+s_0\int_0^\infty \bigg|\frac{\p^2 }{\p t^2}\bigr[{g_1}(\xi,\theta;t+c)\bigr]\bigg|dt.
\end{align*}
Since $|\xi-c|\geq|\xi|$ for all $-\xi_1\leq \xi \leq \xi_2$, we have
\beq \label{g_less_1}
\big|g(\xi,\theta;c)\big|\leq 1.
\eeq
Moreover, we have $\frac{\sinh(\xi-c)}{2a_\ep}=O(1)$, so that it follows
\beq\label{g2}\big| {g_1}(\xi,\theta;c)\big| \leq M, \eeq
for a constant $M$ independent of $\ep$.

In the follows, we show that there is a positive constant $C$ independent of ${\ep},\xi,\theta$ satisfying$$\left|\frac{\p }{\p t}\big[ {g_1}(\xi,\theta;t+c)\big]\Big|_{t=0}\right|,\ \int_0^\infty \bigg|\frac{\p^2}{\p t^2} \big[{g_1}(\xi,\theta;t+c)\big]\bigg|dt
\leq \frac{C}{\sqrt{{\ep}}}.$$
Remind that $w_\theta$ is the function given by $$w_\theta(\xi)=\sqrt{\cosh\xi-\cos\theta}.$$
A straightforward but tedious computation shows
\begin{align*}
\ds&\frac{\p }{\p t} \bigr[{g_1}(\xi,\theta;t+c)\bigr]\Big|_{t=0}=\frac{w_\theta^3(\xi)}{a_\ep}
\frac{3-\cosh^2(c-\xi)-2\cos\theta\cosh(c-\xi)}{4\bigr(\cosh(c-\xi)-\cos\theta\bigr)^{5/2}},\\[2mm]
\ds&\frac{\p^2 }{\p t^2}\bigr[ {g_1}(\xi,\theta;t+c)\bigr]=-\frac{w_\theta^3(\xi)}{a_\ep}
\frac{\sinh\tilde{x}~\bigr(4\cos^2\theta+10\cos\theta\cosh\tilde{x}+\cosh^2\tilde{x}-15\bigr)}{8\bigr(\cosh\tilde{x}-\cos\theta\bigr)^{7/2}}
\end{align*}
with $\tilde{x}=t+c-\xi$.

Thanks to \eqnref{g_less_1}, there exists a constant $C$ independent of $\ep$ such that
\begin{align}
\left|\frac{\p }{\p t} \bigr[{g_1}(\xi,\theta;t+c)\bigr]\Big|_{t=0}\right| & \leq \frac{C}{a_\ep}\left(\left| \frac{3(1-\cosh^2(c-\xi))}{\cosh(c-\xi)-\cos\theta}\right| + 2(\cosh(c-\xi)-\cos\theta) \right) \leq \frac{C}{\sqrt{{\ep}}}.\label{g2prime}
\end{align}
For the inequality of the second derivative of ${g_1}$, let us consider separately the cases $0\leq t \leq1$ and $t>1$.
For $0\leq t \leq1$, we have
\begin{align*}
\ds&\left|\frac{\p^2 }{\p t^2} \bigr[{g_1}(\xi,\theta;t+c)\bigr]\right|\\
&=
\ds\frac{w_\theta(\xi)^3}{a_\ep(\cosh\tilde{x}-\cos\theta)^{5/2}}
\left|\frac{\sinh\tilde{x}(4\cos^2\theta+10\cos\theta\cosh\tilde{x}+\cosh^2\tilde{x}-15)}
{8(\cosh\tilde{x}-\cos\theta)}\right|\\
&\ds\leq \frac{w_\theta(\xi)^3 \sinh \tilde{x}}{a_\ep(\cosh\tilde{x}-\cos\theta)^{5/2}}
\left(
\left| 15\frac{\cosh^2 \tilde{x}-1}{\cosh \tilde{x}-\cos\theta}\right|+\left|4\frac{\cosh^2 \tilde{x}-\cos\theta}{\cosh \tilde{x}-\cos\theta}\right|+10\cosh \tilde{x} \right)\\
&\ds\leq \frac{C}{\sqrt{{\ep}}}\frac{w_\theta(\xi)^3\sinh\tilde{x}}{(\cosh\tilde{x}-\cos\theta)^{5/2}}.
\end{align*}
Hence it follows
\begin{align}
\int_0^1\left| \frac{\p^2 }{\p t^2} \bigr[{g_1}(\xi,\theta;t+c)\bigr] \right| dt &\leq \frac{C}{\sqrt{{\ep}}} w_\theta(\xi)^3 \int_{c-\xi}^{1+c-\xi} \frac{\sinh t}{(\cosh t-\cos\theta)^{5/2}} dt\notag\\
& \leq \frac{C}{\sqrt{{\ep}}}\left( \left| \frac{w_\theta(\xi)^3}{w_\theta(c-\xi)^3}\right|+
\left|\frac{w_\theta(\xi)^3}{w_\theta(1+c-\xi)^3}\right| \right)\label{g2pprime1}
 \leq\frac{C}{\sqrt{{\ep}}}.
\end{align}
For $t>1$, $\left|{\p^2 }/{\p t^2} \bigr[{g_1}(\xi,\theta;t+c)]\right|$ is uniformly bounded by a exponentially decreasing function of $t$, so that we have $$
\int_1^\infty \left|\frac{\p^2 }{\p t^2} \bigr[{g_1}(\xi,\theta;t+c)\bigr]\right| dt\leq C,
$$
for a constant of $C$ independent of $\ep,\xi,\theta$ and $c$.
Using this equation and \eqnref{g2}, \eqnref{g2prime}, and \eqnref{g2pprime1} as well, we prove the lemma.
\qed

\section{Numerical Illustration}
In this section we illustrate the main results with some examples. More precisely, we plot the graphs of $\nabla h$, $\nabla(u-H)$ and their blow-up terms on $\p B_1$. We consider the asymptotic behavior of $\nabla h$ and $\nabla (u-H)$ only on $\p B_1$ because they have the similar behavior on $\p B_2$.
To have precise values of $\nabla h$ and $\nabla (u-H)$, we use their exact solution formulas derived in section \ref{subsec:exactseries}. On the other hand, gradient blow-up terms are simple elementary functions which are easy to compute as explained in section \ref{subsec:gradientblowup}.

\subsection{Exact solution}\label{subsec:exactseries}
The tangential component of $\nabla h$ is zero on $\p B_1$ due to the second condition in \eqnref{h:eqn} and the normal component of $\nabla h$ has the following exact {solution}:

\begin{align}\notag
\ds{\p_\nu h}\big|_{\p B_1}(\theta)&= \frac{\cosh\xi_1-\cos\theta}{a_\ep} \frac{\p}{\p \xi} h(\xi,\theta)\bigg|_{\xi=-\xi_1}\\\notag
\ds&=\frac{\sqrt{2}(\cosh\xi_1-\cos\theta)^{3/2}}{a_\ep} \sum_{n=0}^\infty \Big(n+\frac{1}{2}\Big)
\ds\left(A_n e^{-(n+\frac{1}{2})\xi_1}- B_n e^{(n+\frac{1}{2})\xi_1}\right) P_n(\cos\theta)\\\label{hnu:numerical}
\qquad &+\frac{\sinh(-\xi_1)(\cosh\xi_1-\cos\theta)^{1/2}}{\sqrt{2}a_\ep}
\sum_{n=0}^\infty
\left(A_n e^{-(n+\frac{1}{2})\xi_1}+ B_n e^{(n+\frac{1}{2})\xi_1}\right) P_n(\cos\theta).
\end{align}

The normal derivative of $(u-H)$ for given entire harmonic function $H$ has {a series representation} similar to \eqnref{hnu:numerical}. Especially for the uniform external field, say $H(\Bx)=E_0x_3$, the exact solution for $u$} can be found in many literatures, for example \cite{lekner11}. To state the solution {explicitly}, we define $$
T(c)=\sum_{n=0}^\infty \frac{(2n+1)(e^{(2n+1)c}+1)}{e^{(2n+1)(\xi_1+\xi_2)}-1}\quad\mbox{for }\quad c>0,
$$
and
$$
V_1=- a_\ep \frac{T_2U_1-T_1U_{12}}{U_1 U_2-U_{12}^2},\quad V_2=a_\ep \frac{T_1U_2-T_2U_{12}}{U_1 U_2-U_{12}^2}$$
with $T_j=T(\xi_j)$, $j=1,2$.
Then the solution $u$ to \eqnref{u:eqn} is represented as follows:
$$
(u-H)(\xi,\theta)=\sqrt{2}E_0\sqrt{\cosh\xi-\cos\theta}\sum_{n=0}^\infty
\left(C_n e^{(n+\frac{1}{2})\xi}+ D_n e^{-(n+\frac{1}{2})\xi}\right) P_n(\cos\theta),
$$
where
\begin{align*}
\ds C_n&=\frac{e^{(2n+1)\xi_1} V_2-V_1-E a_\ep (2n+1)(e^{(2n+1)\xi_1}+1) }{e^{(2n+1)(\xi_1+\xi_2)}-1},\\
\ds D_n&=\frac{e^{(2n+1)\xi_2} V_1-V_2+E a_\ep (2n+1)(e^{(2n+1)\xi_2}+1) }{e^{(2n+1)(\xi_1+\xi_2)}-1}.
\end{align*}
Hence, we have similarly to \eqnref{hnu:numerical} that
\begin{align}
\p_\nu(u-H)\Big|_{\p B_1}(\theta)&=E_0\frac{\sqrt{2}(\cosh\xi_1-\cos\theta)^{3/2}}{a_\ep} \sum_{n=0}^\infty \Big(n+\frac{1}{2}\Big)
\left(C_n e^{-(n+\frac{1}{2})\xi_1}- D_n e^{(n+\frac{1}{2})\xi_1}\right) P_n(\cos\theta)\nonumber \\
&\quad +E_0\frac{\sinh(-\xi_1)(\cosh\xi_1-\cos\theta)^{1/2}}{\sqrt{2}a_\ep}
\sum_{n=0}^\infty
\left(C_n e^{-(n+\frac{1}{2})\xi_1}+ D_n e^{(n+\frac{1}{2})\xi_1}\right) P_n(\cos\theta). \label{u_normal_x3}
\end{align}

\subsection{Gradient blow-up terms}\label{subsec:gradientblowup}
The singular function $h$ satisfies from \eqnref{eqn:1013_2} and \eqnref{eqn:nuexi} that
\begin{gather*}\nabla h\big|_{\p B_1}=\mbox{bounded term}+q_{\p B_1}(\Bx)\hskip .5mm\nu(\Bx),\\[1mm]
q_{\p B_1}(\theta)=  q_h(-\xi_1,\theta)\quad\mbox{for }\theta\in[0,\pi].
\end{gather*}
From \eqnref{gjIj} and the definition of $q_h$ in \eqnref{qhint},
we can easily derive
\begin{align}\label{qh:numerical}
\ds q_{\p B_1}(\theta)&=
\frac{\tilde\mu_1}{2a_\ep (\xi_1+\xi_2)}\left( w_\theta^2(\xi_1)+
\frac{w_\theta^3(\xi_1)}{w_\theta(3\xi_1+2\xi_2)}
\right)+\frac{\tilde\mu_2}{2a_\ep(\xi_1+\xi_2)} \left(\frac{2w_\theta^3(\xi_1)}{w_\theta(\xi_1+2\xi_2)}\right).
\end{align}
Recalling $w_\theta(\xi)=\sqrt{\cosh\xi-\cos\theta}$, we see that the function $q_{\p B_1}(\theta)$ consists of elementary functions which can be easily computed numerically. Similarly, the solution $u$ to \eqnref{u:eqn} for {a} given entire harmonic function $H$ satisfies from Proposition \ref{uasymp:bispherical} that
$$\nabla (u-H)\big|_{\p B_1}=\mbox{bounded term}+C_H q_{\p B_1}(\Bx)\hskip .5mm\nu(\Bx).$$
When an uniform field $H(\Bx)=E_0x_3$ is applied, $C_H$ becomes from Theorem \ref{lem:CH}(b) as follows: $$
C_H=E_0\, \mathcal{Q}_1(r_1,r_2)=
4\pi\tilde{r}^{2}E_0 \left[
(\mu_1+\mu_2)({\pi^2}/{6})+\mu_1 \psi_1(\rtwo)+\mu_2 \psi_1(\rone)\right].
$$

\subsection{Examples}

\noindent\textbf{Data Acquisition}
We numerically compute $\p_\nu h$ and $\p_\nu (u-H)$ based on the exact solution \eqnref{hnu:numerical} and \eqnref{u_normal_x3}. It is worth to mention the difficulty in the numerical computation of \eqnref{hnu:numerical} and \eqnref{u_normal_x3}. Since the term $e^{-2n(\xi_1+\xi_2)}$ decays very slowly for small $\ep$, the cost in numerical computation becomes very high. For instance, in Example 1, we evaluate the summation for $n\leq 5\cdot 10^3$ to compute within a relative tolerance $10^{-5}$ when $\ep=5\cdot 10^{-5}$.
On the other hand, the gradient blow-terms $q_{\p B_1}$ and $C_H q_{\p B_1}$ are consists of simple elementary functions, see \eqnref{qh:numerical}. Hence, the computing cost is extremely low.
For all examples, the radii of the two sphere are $r_1=3$ and $r_2=2$.

\smallskip

\begin{figure}[!htp]
\begin{center}
\epsfig{figure=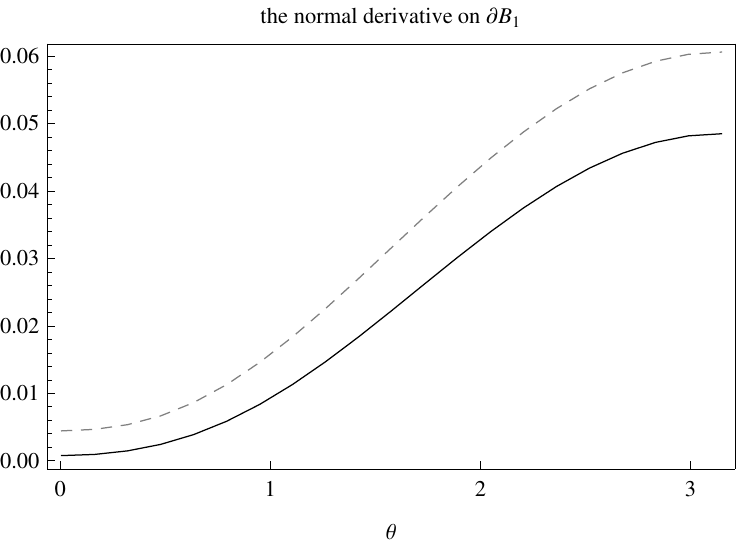,width=4.8cm}\hskip .2cm
\epsfig{figure=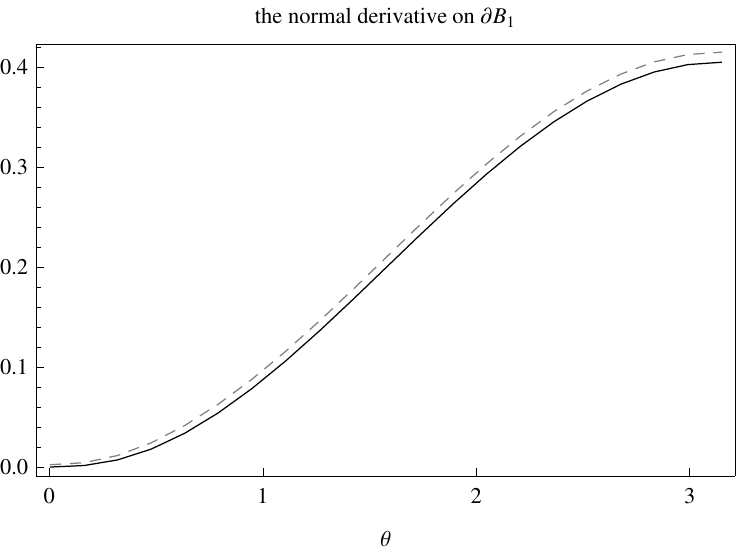,width=4.8cm}\hskip .2cm
\epsfig{figure=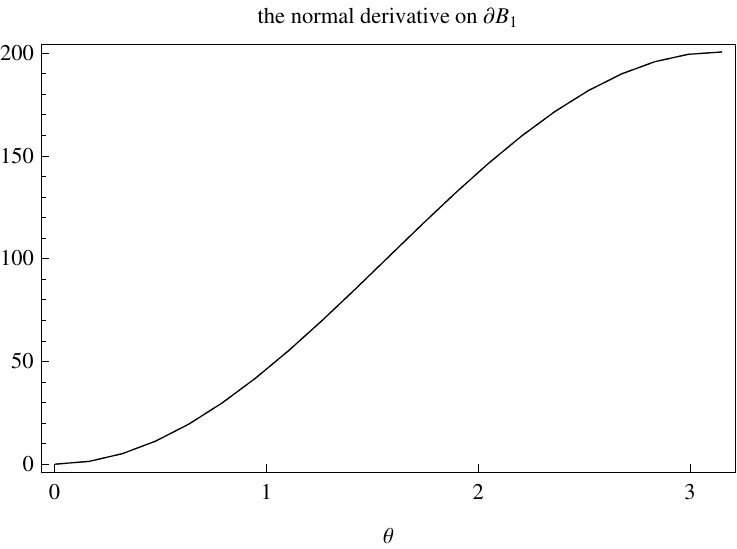,width=4.8cm}\hskip .2cm\\
\epsfig{figure=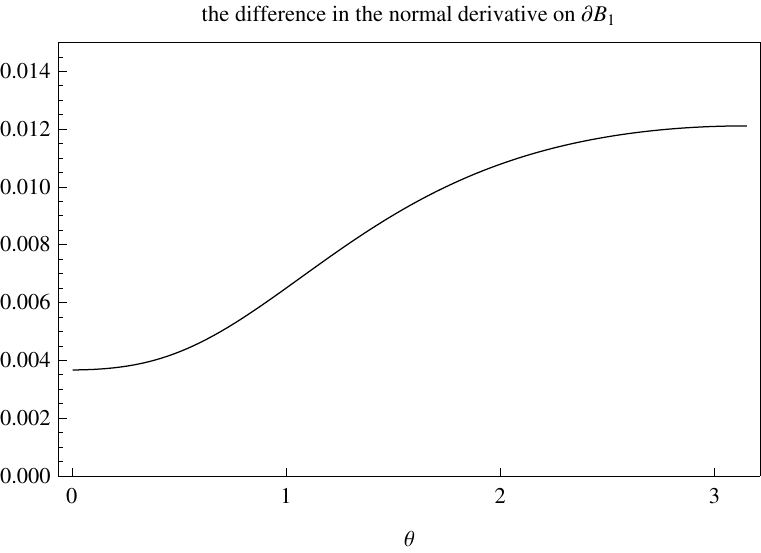,width=4.8cm}\hskip .2cm
\epsfig{figure=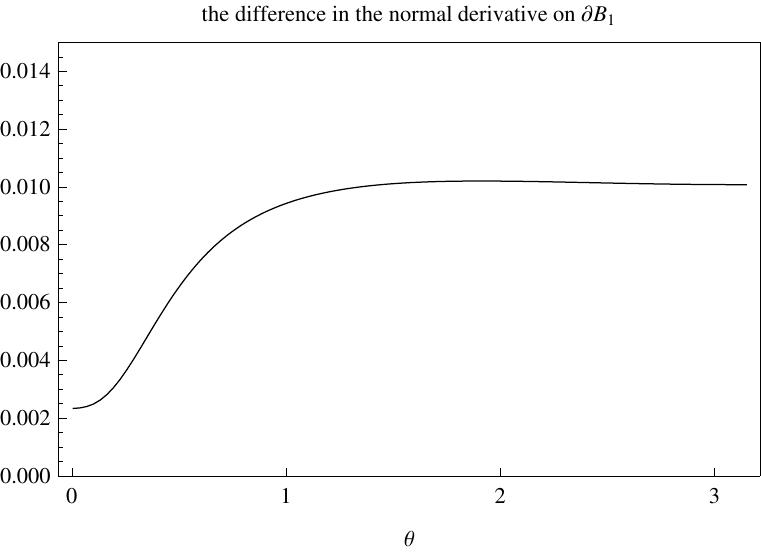,width=4.8cm}\hskip .2cm
\epsfig{figure=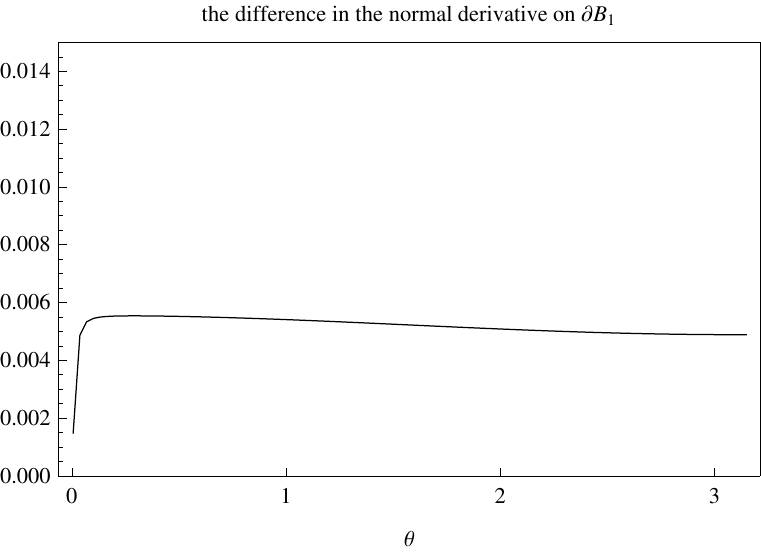,width=4.8cm}\hskip .2cm
\end{center}
\vskip -.8cm
\caption{The graphs of $\p_\nu h|_{\p B_1}$ (dashed), its blow-up term $q_{\p B_1}$ (solid) in the first row, and their difference $|\p_\nu h|_{\p B_1}-q_{\p B_1}|$ in the second row.
The distance $\epsilon$ is $0.5, 0.05, 0.00005$ from left to right columns. }
\label{fig_ep_decr}
\end{figure}

\noindent\textbf{Example 1.}
 In Fig.\;{\rm\ref{fig_ep_decr}}, we compare $\p_\nu h|_{\p B_1}$ and its blow-up term $q_{\p B_1}$ when
 $\epsilon$ takes the values $0.5, 0.05, 0.00005$ from left to right columns. We plot $\p_\nu h|_{\p B_1}$ (dashed graph) and $q_{\p B_1}$ (solid graph) in the first row and the difference between them in the second row. Note that while the range of $y$-axis in the first row becomes huge for small $\ep$, that in the second row is fixed. It means that the magnitudes of both $\p_\nu h|_{\p B_1}$ and $q_{\p B_1}$ increase as $\ep$ decreases, but the difference between them decreases. Hence the blow-up term $q_{\p B_1}$ represents $\p_\nu h|_{\p B_1}$ better when $\ep$ is smaller.

\smallskip
\noindent\textbf{Example 2.}
In Table \ref{table: nor_deri}, we provide the values of $\p_\nu(u-H)|_{\p B_1}$ and its blow-up term $C_H q_{\p B_1}$ for various $\ep$ and various bispherical coordinates values $\theta$ when an uniform external field $H(\Bx)=x_3$ is applied. The difference between $\p_\nu(u-H)|_{\p B_1}$ and its blow-up terms is of almost constant magnitude while the value of $\p_\nu(u-H)|_{\p B_1}$ is huge near $\theta=\pi$ for small $\ep$.

\begin{table}[!p]
\centering 
\renewcommand{\arraystretch}{1}
\setlength{\tabcolsep}{9pt}
\begin{tabular}{|c |c r r | c r r |} 
\hline
${\theta}$ (unit in $\pi$) & \multicolumn{1}{c}{$\ep$}&  \multicolumn{1}{c}{\ ${\p_\nu (u-H)}$} & \multicolumn{1}{c}{$C_H q_{\p B_1}$} &\multicolumn{1}{|c}{ $\ep$\quad}
&\multicolumn{1}{c}{${\p_\nu (u-H)}$} &\multicolumn{1}{c|}{ $C_H q_{\p B_1}$}\\ [0.5ex] 
\hline \hline
$0$ &  1      & $2.118$ &0.103   &   0.00005     & 1.4 & $0.3$ \\
$0.15$ &          & $2.380$ &0.195   &        & 1279.4 & $1278.8$ \\
$0.30$ &        & $2.754$ & 0.505  &         & 4838.5 & $4837.8$\\
$0.45$ &        & $3.200$ & 1.031  &           & 9901.5 & $9900.8$\\
$0.60$ &        & $3.807$ & 1.668  &          & 15365.8 & $15364.1$\\
$0.75$ &        & $4.420$ & 2.251  &           & 20037.4 & $20036.6$\\
$0.90$ &        & $4.827$ & 2.621  &            & 22900.7 & $22900.0$\\
$1.00$ &        & $4.911$ & 2.700  &           & 23475.2 & $23474.5$\\
 & & & & & & \\
$0$    & 0.5    & $1.875$ & 0.089  & 0.000005   & 1.4 & $0.02$\\
$0.15$ &          & $2.298$ &0.278   &        & 10896 & $10896$ \\
$0.30$ &        & $2.877$ & 0.963  &           & 41211 & $41211$\\
$0.45$ &        & $3.938$ & 2.131  &          & 84337 & $84336$\\
$0.60$ &        & $5.360$ & 3.511  &           & 130871 & $130871$\\
$0.75$ &        & $6.688$ & 4.746  &           & 170671 & $170670$\\
$0.90$ &        & $7.529$ & 5.519  &           & 195060 & $195060$\\
$1.00$ &        & $7.700$ & 5.675  &           & 199954 & $199953$\\
 & & & & & & \\
$0$    &  0.05  & $1.586$ & 0.059  & $5\times 10^{-7}$   & 1.4 & $0.02$\\
$0.15$ &          & $3.346$ &2.096   &        & 94900 & $94900$ \\
$0.30$ &        & $10.339$ & 9.117  &           & 358914 & $358914$\\
$0.45$ &        & 20.811 & $19.469$  &          & 734490 & $734490$\\
$0.60$ &        & $32.142$ & 30.718  &           & $1.13976 \times 10^6$ & $1.13976 \times 10^6$\\
$0.75$ &        & $41.838$ & 40.358  &           & $1.48637 \times 10^6$ & $1.48637 \times 10^6$\\
$0.90$ &        & $47.781$ & 46.270  &           & $1.69878 \times 10^6$ & $1.69878 \times 10^6$\\
$1.00$ &        & $48.973$ & 47.456  &           & $1.74140 \times 10^6$ & $1.74139 \times 10^6$\\
 & & & & & & \\
$0$    &  0.005  & 1.489 & $0.043$ & $5\times 10^{-8}$   & 1.4 & $0.02$\\
$0.15$ &         & 19.800 & $18.855$   &        & $840495$ & $840495$ \\
$0.30$ &          & 74.331 & $73.323$&           & $3.17876 \times 10^6$ & $3.17876 \times 10^6$\\
$0.45$ &        & 151.926 & $150.880$  &          & $6.50509 \times 10^6$ & $6.50509 \times 10^6$\\
$0.60$ &        & 235.658 & $234.579$  &           & $1.00944 \times 10^7$ & $1.00944 \times 10^7$\\
$0.75$ &        & 307.272 & $306.166$  &           & $1.31642 \times 10^7$ & $1.31642 \times 10^7$\\
$0.90$ &        & 351.158 & $350.035$ &           & $1.50454 \times 10^7$ & $1.50454 \times 10^7$\\
$1.00$ &        & 359.962 & $358.836$ &           & $1.54228 \times 10^7$ & $1.54228 \times 10^7$\\
 & & & & & & \\
$0$    &  0.0005  & 1.44 & $0.03$  & $5 \times 10^{-9}$   & 1.4 & $0.02$\\
$0.15$ &         & 155.01 & $154.22$   &        & $7.54254\times 10^6$ & $7.54254\times 10^6$ \\
$0.30$ &        & 588.88 & $585.07$  &           & $2.85260\times 10^7$ & $2.85260\times 10^7$\\
$0.45$ &        & 1198.81 & $1197.98$ &          & $5.83762\times 10^7$ & $5.83762\times 10^7$\\
$0.60$ &        & 1860.20 & $1859.35$  &           & $9.05863\times 10^7$ & $9.05863\times 10^7$\\
$0.75$ &        & 2425.87 & $2425.00$  &           & $1.18135 \times 10^8$ & $1.18135\times 10^8$\\
$0.90$ &        & 2772.52 & $2771.63$  &           & $1.35017\times 10^8$ & $1.35017\times 10^8$\\
$1.00$ &        & 2842.07 & $2841.18$ &           & $1.38404\times 10^8$ & $1.38404\times 10^8$\\
\hline 
\end{tabular}\caption{comparison between the exact normal derivative $\p_\nu (u-H)|_{\p B_1}$ and its blow-up term $C_Hq_{\p B_1}$} \label{table: nor_deri}
\end{table}


\section{Conclusion}
In this paper we provided an asymptotic analysis for the superfocusing of the electric field due to the presence of two nearly touching perfectly conducting spheres. We expressed explicitly and completely the blow-up term of the electric field with the rigorous proof. The main ideas of this paper come from, firstly, the solution by separation of variables in the bispherical coordinates and, secondly, the idea to approximate the series solution by an integral function using the Euler-Maclaurin formula and, thirdly, the recent decomposition method to separate the blow-up term and the regular term in the electric field. The derived asymptotic formula is valid in the whole exterior region of the two spheres, and it explicitly characterizes superfocusing of the electric field.

\end{document}